\documentclass{amsart}

\usepackage{amssymb}
\usepackage{amsmath}
\usepackage{amsthm}
\usepackage{amsfonts}

\usepackage{a4wide}
\usepackage{longtable}
\usepackage{multirow}
\usepackage{setspace}

\newtheorem{theorem}{Theorem}[section]

\theoremstyle{definition}

\newtheorem{example}{Example}[section]

\numberwithin{equation}{section}

\begin{document}


\title[Pattern recognition via Artin transfers and finite 2-class towers]
{The strategy of pattern recognition via Artin transfers \\
applied to finite towers of 2-class fields}

\author{Daniel C. Mayer}
\address{Naglergasse 53\\8010 Graz\\Austria}
\email{algebraic.number.theory@algebra.at}
\urladdr{http://www.algebra.at}

\thanks{Research supported by the Austrian Science Fund (FWF): project P 26008-N25}


\date{June 19, 2019}


\begin{abstract}
The isomorphism type of the Galois group of the \(2\)-class field tower
of quadratic number fields having a \(2\)-class group with abelian type invariants \((4,4)\)
is determined by means of information on the transfer of \(2\)-classes to unramified abelian \(2\)-extensions,
collected in the Artin pattern.
In recent investigations by Benjamin and Snyder,
the length of the tower of such fields has turned out to be dependent on
the rank of the \(2\)-class group of the first Hilbert \(2\)-class field.
Significant progress is achieved by extending the pool of possible metabelian \(2\)-groups
of the second Hilbert \(2\)-class field
from the SmallGroups database, resp. Hall-Senior classification,
with the aid of the \(p\)-group generation algorithm,
and sifting the pool by means of pattern recognition.
\end{abstract}

\maketitle


\section{Introduction}
\label{s:Intro}
\noindent
Let \(K\) be a quadratic number field
with \(2\)-class group \(\mathrm{Cl}_2(K)\simeq C_4\times C_4\),
that is, with abelian type invariants \((4,4)\).
The strategy of \textit{pattern recognition via Artin transfers}
consists of three steps.

First, the number theoretic Artin pattern \(\mathrm{AP}(K)\) of \(K\) is defined as
the collection \((\varkappa(K),\tau(K))\)
of all kernels \(\varkappa(K)=(\ker(T_{K,N}))_N\) 
and targets \(\tau(K)=(\mathrm{Cl}_2(N))_N\) 
of \textit{class extension} homomorphisms
\(T_{K,N}:\,\mathrm{Cl}_2(K)\to\mathrm{Cl}_2(N)\)
from \(K\) to unramified abelian extensions \(K<N<\mathrm{F}_2^1(K)\)
of \(K\) within the Hilbert \(2\)-class field \(\mathrm{F}_2^1(K)\) of \(K\)
\cite{Ma2011}.
It is determined by means of the
class field theoretic routines of the computational algebra system Magma
\cite{BCP,BCFS,MAGMA}.

Second, the number theoretic Artin pattern \(\mathrm{AP}(K)\) of \(K\)
is interpreted as the group theoretic Artin pattern \(\mathrm{AP}(G)\) of the
metabelian Galois group \(G=\mathrm{Gal}(\mathrm{F}_2^2(K)/K)\)
of the second Hilbert \(2\)-class field \(\mathrm{F}_2^2(K)\) of \(K\),
which consists of all kernels \(\varkappa(G)=(\ker(T_{G,S}))_S\)
and targets \(\tau(G)=(S/S^\prime)_S\)
of \textit{Artin transfer} homomorphisms
\(T_{G,S}:\,G/G^\prime\to S/S^\prime\)
from \(G\) to normal subgroups \(G^\prime<S<G\)
containing the commutator subgroup \(G^\prime\) of \(G\)
\cite{Ma2010},
which is isomorphic to the \(2\)-class group
\(\mathrm{Cl}_2(\mathrm{F}_2^1(K))\simeq\mathrm{Gal}(\mathrm{F}_2^2(K)/\mathrm{F}_2^1(K))\)
of the Hilbert \(2\)-class field \(\mathrm{F}_2^1(K)\) of \(K\).
If the unramified abelian extension \(N/K\) is the fixed field \(N=\mathrm{Fix}(S)\)
of the subgroup \(S<G\),
in the sense of the Galois correspondence,
then \(S/S^\prime\simeq\mathrm{Gal}(\mathrm{F}_2^2(K)/N)/\mathrm{Gal}(\mathrm{F}_2^2(K)/\mathrm{F}_2^1(N))\simeq\mathrm{Cl}_2(N)\),
according to the Artin reciprocity law
\cite{Ar1927},
and \(\ker(T_{G,S})\simeq\ker(T_{K,N})\),
by the theory of the Artin transfer
\cite{Ar1929}.

Third, the group theoretic Artin pattern \(\mathrm{AP}(G)\)
is used as a \textit{search pattern} in some pool of metabelian \(2\)-groups \(G\)
with abelianization \(G/G^\prime\simeq C_4\times C_4\),
for instance the Hall-Senior classification
or the SmallGroups library
\cite{BEO},
which contains all finite \(2\)-groups up to order \(2^9\),
or an \textit{extended pool}, constructed selectively with the aid of the
\(p\)-group generation algorithm by Newman
\cite{Nm}
and O'Brien
\cite{Ob}.
Frequently, this \textit{sifting process} yields a unique candidate
for the second \(2\)-class group \(G=\mathrm{Gal}(\mathrm{F}_2^2(K)/K)\)
of the quadratic field \(K\)
\cite{Ma2012},
and the relation rank \(d_2(G)\) of \(G\) admits a decision about
the length \(\ell_2(K)\) of the \(2\)-class field tower of \(K\),
according to the \textit{Shafarevich cohomology criterion}
\cite[\S\S\ 1.2--1.3, pp. 75--76]{Ma2016b}
in \S\
\ref{s:Shafarevich}.

Concerning the layout of this paper:
The lattice of normal subgroups of \(G/G^\prime\simeq C_4\times C_4\)
is given as a foundation in \S\
\ref{s:NormalLattice},
imaginary base fields are classified in \S\
\ref{s:Imaginary44},
 real base fields in \S\
\ref{s:Real44}.

An extensive class of imaginary quadratic fields,
presented in Theorems
\ref{thm:GroundState},
\ref{thm:GroundStateVar} and
\ref{thm:ExcitedState},
gives the opportunity to endow suitable subsets
of the amorphous pool of finite \(2\)-groups
with the structure of \textit{descendant trees} by parent-child relations
\cite{Ma2016a}.
The periodicity of the branches of such trees with fixed coclass
admits the characterization of infinitely many groups
by a \textit{parametrized system of invariants}
in \S\
\ref{ss:FirstTree},
and the utilization of the fundamental \textit{principles} of
\textit{monotony, stabilization, polarization} and \textit{mainline}
in the proof of many theorems.


\section{Normal lattice of \(G/G^\prime\simeq C_4\times C_4\)}
\label{s:NormalLattice}
\noindent
Let \(G\) be a finite \(2\)-group with two generators \(x,y\) such that
\(G=\langle x,y\rangle\) and \(G/G^\prime\simeq C_4\times C_4\),
i.e., \(x^4,y^4\in G^\prime\),
where \(G^\prime\) denotes the \textit{derived subgroup} of \(G\).
Let \(\tau_0=(G/G^\prime)\) (0th layer).
In the following description of the \textit{normal lattice} of \(G\),
we use a notation similar to
\cite[\S\ 2]{BeSn}.

The \textit{Frattini subgroup} \(J_0:=\Phi(G)=\langle x^2,y^2,G^\prime\rangle\) of \(G\)
is the intersection of the three \textit{maximal subgroups}
\[
H_1:=\langle x,J_0\rangle,\quad H_2:=\langle xy,J_0\rangle,\quad H_3:=\langle y,J_0\rangle,\quad
\text{ with 1st layer Artin pattern }(\varkappa_1,\tau_1).
\]
They are of index \(2\) in \(G\) and can also be viewed as
\[
H_1=\langle x,y^2,G^\prime\rangle,\quad H_2=\langle xy,y^2,G^\prime\rangle,\quad H_3=\langle y,x^2,G^\prime\rangle,
\]
which shows that \(H_i/G^\prime\simeq C_4\times C_2\), for each \(1\le i\le 3\),
whereas \(J_0/G^\prime\simeq C_2\times C_2\).

In addition to \(J_0\), there are six subgroups \(J_{ik}\)
with \(1\le i\le 3\), \(1\le k\le 2\), of index \(4\) in \(G\): \\
two subgroups of \(H_1\),
\(\quad J_{11}:=\langle x,G^\prime\rangle,\quad J_{12}:=\langle xy^2,G^\prime\rangle\), \\
two subgroups of \(H_2\),
\(\quad J_{21}:=\langle xy,G^\prime\rangle,\quad J_{22}:=\langle xy^3,G^\prime\rangle\), \quad and \\
two subgroups of \(H_3\),
\(\quad J_{31}:=\langle y,G^\prime\rangle,\quad J_{32}:=\langle x^2y,G^\prime\rangle,\quad
\text{ with 2nd layer Artin pattern }(\varkappa_2,\tau_2)\). \\
They all have a cyclic quotient \(J_{ik}/G^\prime\simeq C_4\).

Finally, there are three subgroups of index \(8\) in \(G\):
\[
K_1:=\langle x^2,G^\prime\rangle,\quad K_2:=\langle x^2y^2,G^\prime\rangle,\quad K_3:=\langle y^2,G^\prime\rangle,\quad
\text{ with 3rd layer Artin pattern }(\varkappa_3,\tau_3).
\]
They have quotient \(K_i/G^\prime\simeq C_2\), for each \(1\le i\le 3\).
Additionally, let \(\tau_4=(G^\prime/G^{\prime\prime})\) (4th layer).

Logarithmic abelian type invariants are used for the components of \(\tau\),
for instance,
\((4111)\) or even \((41^3)\), with formal exponents indicating iteration, instead of \((16,2,2,2)\).

Among the remaining members of the normal lattice of \(G\), we mention
the terms \(\gamma_j\) of the \textit{lower central series} with
\(\gamma_1:=G\), and \(\gamma_j:=\lbrack\gamma_{j-1},G\rbrack\), for \(2\le j\le c\),
and the terms \(\zeta_j\) of the \textit{upper central series} with 
\(\zeta_0:=1\), and \(\zeta_j/\zeta_{j-1}:=\mathrm{center}(G/\zeta_{j-1})\), for \(1\le j\le c-1\).
In both cases, \(c=\mathrm{cl}(G)\) denotes the \textit{nilpotency class} of \(G\)
such that \(\gamma_{c+1}=1\) and \(\zeta_c=G\).
If \(G\) is of order \(2^n\), then \(r=\mathrm{cc}(G)=n-c\) is the \textit{coclass} of \(G\).


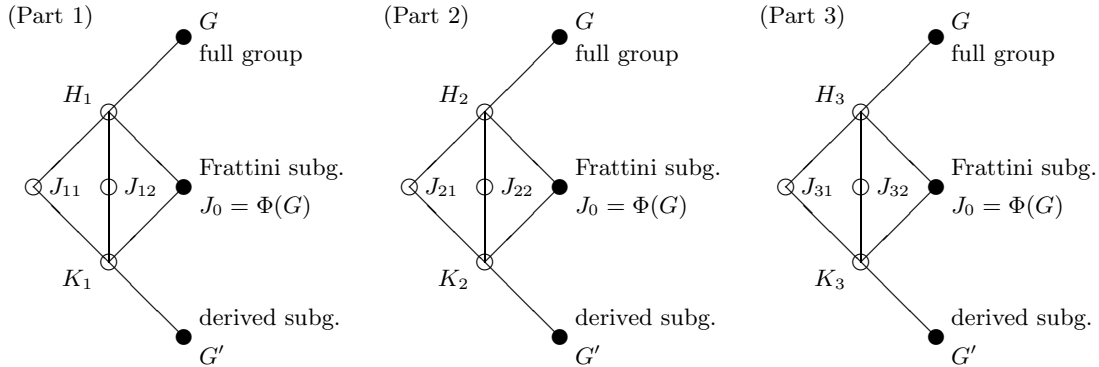
\begin{figure}[ht]
\caption{Normal lattice of \(G/G^\prime\simeq C_4\times C_4\)}
\label{fig:NormalLattice}

{\small

\setlength{\unitlength}{1.0cm}
\begin{picture}(10,6)(-8,-10)



\put(-9.2,-4.9){\makebox(0,0)[rb]{(Part 1)}}
\put(-8,-5){\circle*{0.2}}
\put(-7.8,-4.9){\makebox(0,0)[lb]{\(G\)}}
\put(-7.8,-5.1){\makebox(0,0)[lt]{full group}}

\put(-8,-5){\line(-1,-1){2}}

\put(-9,-6){\circle{0.2}}
\put(-9.2,-5.9){\makebox(0,0)[rb]{\(H_1\)}}

\put(-9,-6){\line(0,-1){2}}

\put(-9,-6){\line(1,-1){1}}

\put(-10,-7){\circle{0.2}}
\put(-9.8,-7){\makebox(0,0)[lc]{\(J_{11}\)}}
\put(-9,-7){\circle{0.2}}
\put(-8.8,-7){\makebox(0,0)[lc]{\(J_{12}\)}}

\put(-7.8,-6.9){\makebox(0,0)[lb]{Frattini subg.}}
\put(-8,-7){\circle*{0.2}}
\put(-7.8,-7.1){\makebox(0,0)[lt]{\(J_0=\Phi(G)\)}}

\put(-8,-7){\line(-1,-1){1}}

\put(-9,-8){\circle{0.2}}
\put(-9.2,-8.1){\makebox(0,0)[rt]{\(K_1\)}}

\put(-8,-9){\line(-1,1){2}}

\put(-8,-9){\circle*{0.2}}
\put(-7.8,-8.9){\makebox(0,0)[lb]{derived subg.}}
\put(-7.8,-9.1){\makebox(0,0)[lt]{\(G^\prime\)}}



\put(-4.2,-4.9){\makebox(0,0)[rb]{(Part 2)}}
\put(-3,-5){\circle*{0.2}}
\put(-2.8,-4.9){\makebox(0,0)[lb]{\(G\)}}
\put(-2.8,-5.1){\makebox(0,0)[lt]{full group}}

\put(-3,-5){\line(-1,-1){2}}

\put(-4,-6){\circle{0.2}}
\put(-4.2,-5.9){\makebox(0,0)[rb]{\(H_2\)}}

\put(-4,-6){\line(0,-1){2}}

\put(-4,-6){\line(1,-1){1}}

\put(-5,-7){\circle{0.2}}
\put(-4.8,-7){\makebox(0,0)[lc]{\(J_{21}\)}}
\put(-4,-7){\circle{0.2}}
\put(-3.8,-7){\makebox(0,0)[lc]{\(J_{22}\)}}

\put(-2.8,-6.9){\makebox(0,0)[lb]{Frattini subg.}}
\put(-3,-7){\circle*{0.2}}
\put(-2.8,-7.1){\makebox(0,0)[lt]{\(J_0=\Phi(G)\)}}

\put(-3,-7){\line(-1,-1){1}}

\put(-4,-8){\circle{0.2}}
\put(-4.2,-8.1){\makebox(0,0)[rt]{\(K_2\)}}

\put(-3,-9){\line(-1,1){2}}

\put(-3,-9){\circle*{0.2}}
\put(-2.8,-8.9){\makebox(0,0)[lb]{derived subg.}}
\put(-2.8,-9.1){\makebox(0,0)[lt]{\(G^\prime\)}}



\put(0.8,-4.9){\makebox(0,0)[rb]{(Part 3)}}
\put(2,-5){\circle*{0.2}}
\put(2.2,-4.9){\makebox(0,0)[lb]{\(G\)}}
\put(2.2,-5.1){\makebox(0,0)[lt]{full group}}

\put(2,-5){\line(-1,-1){2}}

\put(1,-6){\circle{0.2}}
\put(0.8,-5.9){\makebox(0,0)[rb]{\(H_3\)}}

\put(1,-6){\line(0,-1){2}}

\put(1,-6){\line(1,-1){1}}

\put(0,-7){\circle{0.2}}
\put(0.2,-7){\makebox(0,0)[lc]{\(J_{31}\)}}
\put(1,-7){\circle{0.2}}
\put(1.2,-7){\makebox(0,0)[lc]{\(J_{32}\)}}

\put(2.2,-6.9){\makebox(0,0)[lb]{Frattini subg.}}
\put(2,-7){\circle*{0.2}}
\put(2.2,-7.1){\makebox(0,0)[lt]{\(J_0=\Phi(G)\)}}

\put(2,-7){\line(-1,-1){1}}

\put(1,-8){\circle{0.2}}
\put(0.8,-8.1){\makebox(0,0)[rt]{\(K_3\)}}

\put(2,-9){\line(-1,1){2}}

\put(2,-9){\circle*{0.2}}
\put(2.2,-8.9){\makebox(0,0)[lb]{derived subg.}}
\put(2.2,-9.1){\makebox(0,0)[lt]{\(G^\prime\)}}


\end{picture}

}

\end{figure}

Figure
\ref{fig:NormalLattice}
shows the three parts of the normal lattice of \(G\) above the derived subgroup \(G^\prime\).


\section{Shafarevich cohomology criterion}
\label{s:Shafarevich}
\noindent
Let \(p\) be a prime number and \(G\) be a pro-\(p\) group.
Then the finite field \(\mathbb{F}_p\) is a trivial \(G\)-module,
and two crucial cohomological invariants of \(G\) are
the \textit{generator rank} \(d_1:=d_{p,1}(G):=\dim_{\mathbb{F}_p}H^1(G,\mathbb{F}_p)\) and
the \textit{relation rank} \(d_2:=d_{p,2}(G):=\dim_{\mathbb{F}_p}H^2(G,\mathbb{F}_p)\).
Shafarevich has established a necessary criterion for the Galois group
\(G=\mathrm{Gal}(\mathrm{F}_p^\infty(K)/K)\)
of the Hilbert \(p\)-class field tower \(\mathrm{F}_p^\infty(K)\) of a number field \(K\)
in terms of the generator rank and relation rank of \(G\)
\cite[Thm. 5.1, p. 28]{Ma2015}:
\[
d_1\le d_2\le d_1+r+\theta,
\]
where \(r\) and \(\theta\) are defined in the following way:
if \((r_1,r_2)\) denotes the signature of \(K\),
then the torsion free \textit{Dirichlet unit rank} of \(K\) is given by \(r=r_1+r_2-1\)
and
\[
\theta=
\begin{cases}
1 & \text{if } K \text{ contains a primitive } p\text{-th root of unity}, \\
0 & \text{otherwise}.
\end{cases}
\]
Applied to the particular case of \(p=2\)
and the \(2\)-class field tower group \(G=\mathrm{Gal}(\mathrm{F}_2^\infty(K)/K)\)
of a quadratic number field \(K\) with \(2\)-class group
\(G/G^\prime\simeq\mathrm{Cl}_2(K)\simeq C_4\times C_4\),
and thus with Frattini quotient \(G/\Phi(G)\simeq C_2\times C_2\),
the Shafarevich criterion becomes
\[
2\le d_2\le 3+r,
\]
since the generator rank is \(d_1=2\), according to the Burnside basis theorem,
\(\theta=1\) trivially, and the Dirichlet unit rank is
\[
r=
\begin{cases}
1 & \text{for real } K \text{ with } (r_1,r_2)=(2,0), \\
0 & \text{for imaginary } K \text{ with } (r_1,r_2)=(0,1).
\end{cases}
\]


\section{Imaginary quadratic fields of type \((4,4)\)}
\label{s:Imaginary44}
\noindent
According to Benjamin and Snyder
\cite[\S\ 5, p. 1183]{BeSn},
the discriminant of an imaginary quadratic field \(K\)
with \(2\)-class group \(\mathrm{Cl}_2(K)\) of type \((4,4)\)
is the product \(d_K=d_1d_2d_3\) of three distinct prime discriminants
\(d_1<0\) and \(d_2,d_3>0\),
and the three unramified quadratic extensions \(N_i:=\mathrm{Fix}(H_i)\) of \(K\)
are given by \(N_i=K(\sqrt{d_i})\), for \(1\le i\le 3\),
such that generally \(\#\mathrm{Cl}_2(N_i)>16=\#\mathrm{Cl}_2(K)\) for \(i\in\lbrace 2,3\rbrace\)
\cite[Lem. 5.1]{BeSn}.
The remaining extension \(N_1\) is responsible
for the \(2\)-rank of \(\tau_4=\mathrm{Cl}_2(\mathrm{F}_2^1(K))\).
In \S\
\ref{ss:Length2},
we first consider the case that \(\#\mathrm{Cl}_2(N_1)=16\)
and \(\#\mathrm{Cl}_2(N_{1k})=16\)
for the two unramified cyclic quartic extensions
\(N_{1k}:=\mathrm{Fix}(J_{1k})\), \(k\in\lbrace 1,2\rbrace\), of \(K\)
which contain \(N_1\),
whence \(\tau_4=\mathrm{Cl}_2(\mathrm{F}_2^1(K))\) has rank \(2\) by
\cite[Cor. 3.2, p. 1182]{BeSn}.


\subsection{Metabelian towers of Hilbert \(2\)-class fields}
\label{ss:Length2}
\noindent
In Table
\ref{tbl:Imaginary2},
some discriminants \(0>d_K>-10^5\) of imaginary quadratic fields \(K\)
with \(\mathrm{Cl}_2(K)\simeq C_4\times C_4\) are shown together with
the factorization of \(d_K\) into three prime discriminants
and the Galois group \(\mathrm{Gal}(\mathrm{F}_2^\infty(K)/K)\)
of the metabelian \(2\)-class field tower of \(K\)
in the notation of the SmallGroups library
\cite{BEO}
and ANUPQ package
\cite{GNO}.


\renewcommand{\arraystretch}{1.1}

\begin{table}[ht]
\caption{Discriminants of imaginary quadratic fields \(K\) with \(\ell_2(K)=2\)}
\label{tbl:Imaginary2}
\begin{center}
\begin{tabular}{|r|c||c|c|}
\hline
 \(d_K\)      & Factorization             & \(\mathrm{Gal}(\mathrm{F}_2^\infty(K)/K)\) & \(\tau_4\) \\
\hline
  \(-6\,052\) & \((-4)\cdot 17\cdot 89\)  & \(\langle 128,31\rangle\)                  & \((21)\)   \\
  \(-8\,103\) & \((-3)\cdot 37\cdot 73\)  & \(\langle 256,285\rangle\)                 & \((31)\)   \\
 \(-11\,972\) & \((-4)\cdot 41\cdot 73\)  & \(\langle 256,286\rangle\)                 & \((31)\)   \\
 \(-14\,547\) & \((-3)\cdot 13\cdot 373\) & \(\langle 512,1605\rangle-\#1;2\)          & \((51)\)   \\
 \(-17\,476\) & \((-4)\cdot 17\cdot 257\) & \(\langle 512,1598\rangle\)                & \((41)\)   \\
\hline
\end{tabular}
\end{center}
\end{table}

\noindent
The root paths of the occurring metabelian \(2\)-groups with orders up to \(1024\) are given by
\eqref{eqn:RootPathI2}.

\begin{equation}
\label{eqn:RootPathI2}
\begin{aligned}
\langle 128,31\rangle\stackrel{2}{\to}\langle 32,2\rangle\to\langle 16,2\rangle \\
\langle 256,285 \text{ or } 286\rangle\to\langle 128,35\rangle\stackrel{2}{\to}\langle 32,2\rangle\to\langle 16,2\rangle \\
\#1;2\to\langle 512,1605\rangle\to\langle 256,284\rangle\to\langle 128,35\rangle\stackrel{2}{\to}\langle 32,2\rangle\to\langle 16,2\rangle \\
\langle 512,1598\rangle\to\langle 256,276\rangle\to\langle 128,33\rangle\stackrel{2}{\to}\langle 32,2\rangle\to\langle 16,2\rangle
\end{aligned}
\end{equation}

\noindent
In accordance with the Shafarevich criterion,
all these groups have relation rank \(d_2=3\).


\begin{theorem}
\label{thm:SporadicI} (Sporadic situation.)
Let \(K\) be an imaginary quadratic field
with Artin pattern
\begin{equation}
\label{eqn:APSporadicI}
\begin{aligned}
\tau_0&=(22), \quad \tau_1=(211,221,311), \quad \varkappa_1=(J_0,K_2,K_2), \\
\tau_2&=(22,22,31,31,211,211;221), \quad \varkappa_2=(H_2,H_2,H_2,H_2,H_3,H_3;H_2), \quad \tau_4=(21).
\end{aligned}
\end{equation}
Then the \(2\)-class field tower of \(K\) has length \(\ell_2(K)=2\)
and metabelian Galois group \(\mathrm{Gal}(\mathrm{F}_2^\infty(K)/K)\)
isomorphic to the sporadic group \(\langle 128,31\rangle\) with coclass \(4\).
\end{theorem}


\begin{example}
\label{exm:SporadicI}
The information in Formula
\eqref{eqn:APSporadicI}
occurs rather sparsely as Artin pattern
of imaginary quadratic fields \(K\) with \(2\)-class group of type \((4,4)\).
The absolutely smallest two discriminants \(d_K\) with this pattern are
\begin{equation}
\label{eqn:ExmSporadicI}
\begin{aligned}
-6\,052&=(-4)\cdot 17\cdot 89 \quad 
\text{\cite[Exm. 7.1, pp. 1191--1192]{BeSn}}, \\
-29\,444&=(-4)\cdot 17\cdot 433. \\
\end{aligned}
\end{equation}
\end{example}


\begin{proof}
(Proof of Theorem
\ref{thm:SporadicI})
Coarse selection using \(\tau_0=(22)\) and the Artin pattern \(\tau_1=(211,221,311)\) of the first layer only
yields four candidates of order \(128=2^7\) in
\cite{BEO}.
See Table
\ref{tbl:Comparison}.

However, the unique candidate with correct invariants \(\varkappa_1\), \(\tau_2\), \(\varkappa_2\)
is \(\langle 128,31\rangle\).
All candidates share the common invariant \(\tau_4=(21)\).
So the second layer enables the fine selection.

Candidates of bigger orders \(256=2^8\) and \(512=2^9\)
possess too big components of \(\tau_2\)
and the search can be stopped, according to the \textit{monotony principle}
\cite[\S\ 1.7, Thm. 1.21, p. 79]{Ma2016b}.

Here we have \(\mathrm{Cl}_2(N_1)=\mathrm{Cl}_2(N_{1k})=(211)\), \(k\in\lbrace 1,2\rbrace\),
and indeed \(\tau_4=(21)\)
\cite[Cor. 3.2]{BeSn}.
\end{proof}


\renewcommand{\arraystretch}{1.1}

\begin{table}[ht]
\caption{Comparison of the candidates for Artin pattern \eqref{eqn:APSporadicI}}
\label{tbl:Comparison}
\begin{center}
{\small
\begin{tabular}{|c||c|c|c|c|c|}
\hline
 Vertex                    & \(\tau_1\)      & \(\varkappa_1\) & \(\tau_2\)                   & \(\varkappa_2\)                    & \(\tau_4\) \\
\hline
 \(\langle 128,26\rangle\) & \(211,221,311\) & \(J_0,K_3,J_0\) & \(22,22,31,31,211,211;2111\) & \(H_3,H_3,G,G,G,G;H_3\)            & \((21)\)   \\
 \(\langle 128,31\rangle\) & \(211,221,311\) & \(J_0,K_2,K_2\) & \(22,22,31,31,211,211;221\)  & \(H_2,H_2,H_2,H_2,H_3,H_3;H_2\)    & \((21)\)   \\
 \(\langle 128,34\rangle\) & \(211,221,311\) & \(J_0,K_3,K_2\) & \(22,22,31,31,31,31;221\)    & \(H_3,H_3,H_2,H_2,H_2,H_2;J_0\)    & \((21)\)   \\
 \(\langle 128,35\rangle\) & \(211,221,311\) & \(J_0,K_1,J_0\) & \(22,22,31,31,31,31;221\)    & \(H_3,H_3,H_2,H_2,H_1,H_1;H_1\)    & \((21)\)   \\
\hline
\end{tabular}
}
\end{center}
\end{table}


\begin{theorem}
\label{thm:GroundState} (Ground state.)
Let \(K\) be an imaginary quadratic field
with Artin pattern
\begin{equation}
\label{eqn:APGS}
\begin{aligned}
\tau_0&=(22), \quad \tau_1=(211,221,411), \quad \varkappa_1=(J_0,K_1,K_3), \\
\tau_2&=(22,22,31,31,41,41;321), \quad \varkappa_2=(H_1,H_1,H_2,H_2,H_3,H_3;J_0), \quad \tau_4=(31).
\end{aligned}
\end{equation}
Then the \(2\)-class field tower of \(K\) has length \(\ell_2(K)=2\)
and metabelian Galois group \(\mathrm{Gal}(\mathrm{F}_2^\infty(K)/K)\)
isomorphic to the periodic group \(\langle 256,285\rangle\) (see \S\
\ref{ss:FirstTree})
with coclass \(4\).
\end{theorem}


\begin{example}
\label{exm:GroundState}
The information in Formula
\eqref{eqn:APGS}
seems to be the most frequent Artin pattern
of imaginary quadratic fields \(K\) with \(2\)-class group of type \((4,4)\).
The absolutely smallest ten discriminants \(d_K\) with this pattern are
\begin{equation}
\label{eqn:ExmGS}
\begin{aligned}
-8\,103&=(-3)\cdot 37\cdot 73, \qquad
&-10\,808&=8\cdot (-7)\cdot 193, \\
-12\,104&=(-8)\cdot 17\cdot 89, \qquad
&-12\,207&=(-3)\cdot 13\cdot 313, \\
-13\,727&=(-7)\cdot 37\cdot 53, \qquad
&-14\,155&=5\cdot (-19)\cdot 149, \\
-17\,751&=(-3)\cdot 61\cdot 97, \qquad
&-23\,439&=(-3)\cdot 13\cdot 601, \\
-28\,067&=13\cdot 17\cdot (-127), \qquad
&-33\,464&=8\cdot (-47)\cdot 89.
\end{aligned}
\end{equation}
\end{example}


\begin{proof}
(Proof of Theorem
\ref{thm:GroundState})
Searching the SmallGroups database
\cite{BEO}
for the pattern \(\tau_0=(22)\), \(\tau_1=(211,221,411)\) and \(\tau_2=(22,22,31,31,41,41;321)\)
yields \(4\) groups \(\langle 256,i\rangle\) with identifiers \(284\le i\le 287\).
Only \(\langle 256,285\rangle\) possesses correct capitulation type
\(\varkappa_1=(J_0,K_1,K_3)\) and \(\varkappa_2=(H_1,H_1,H_2,H_2,H_3,H_3;J_0)\)
(Table
\ref{tbl:Variants}).
Here, \(\mathrm{Cl}_2(N_1)=(211)\), \(\mathrm{Cl}_2(N_{1k})=(31)\),
\(\tau_4=(31)\)
\cite[Cor. 3.2]{BeSn}
\end{proof}


\begin{theorem}
\label{thm:GroundStateVar} (Ground state variant.)
Let \(K\) be an imaginary quadratic field
with Artin pattern
\begin{equation}
\label{eqn:APGSVar}
\begin{aligned}
\tau_0&=(22), \quad \tau_1=(211,221,411), \quad \varkappa_1=(J_0,K_1,K_1), \\
\tau_2&=(22,22,31,31,41,41;321), \quad \varkappa_2=(H_1,H_1,H_2,H_2,H_1,H_1;H_1), \quad \tau_4=(31).
\end{aligned}
\end{equation}
Then the \(2\)-class field tower of \(K\) has length \(\ell_2(K)=2\)
and metabelian Galois group \(\mathrm{Gal}(\mathrm{F}_2^\infty(K)/K)\)
isomorphic to the periodic group \(\langle 256,286\rangle\) (see \S\
\ref{ss:FirstTree})
with coclass \(4\).
\end{theorem}


\begin{example}
\label{exm:GroundStateVar}
Only a single imaginary quadratic fields \(K\) with \(2\)-class group of type \((4,4)\)
and Artin pattern in Formula
\eqref{eqn:APGSVar}
is known up to now.
The discriminant \(d_K\) of this field is
\(-11\,972=(-4)\cdot 41\cdot 73\).
\end{example}


\begin{proof}
(Proof of Theorem
\ref{thm:GroundStateVar})
As in the previous proof,
searching the SmallGroups database
\cite{BEO}
for the pattern \(\tau_0=(22)\), \(\tau_1=(211,221,411)\) and \(\tau_2=(22,22,31,31,41,41;321)\)
yields \(4\) groups \(\langle 256,i\rangle\) with \(284\le i\le 287\).
However, only \(\langle 256,286\rangle\) possesses the correct capitulation type
\(\varkappa_1=(J_0,K_1,K_1)\) with two equal kernels and \(\varkappa_2=(H_1,H_1,H_2,H_2,H_1,H_1;H_1)\).
See Table
\ref{tbl:Variants}.
\end{proof}


\begin{theorem}
\label{thm:ExcitedState} (Excited state.)
Let \(K\) be an imaginary quadratic field
with Artin pattern
\begin{equation}
\label{eqn:APES}
\begin{aligned}
\tau_0&=(22), \quad \tau_1=(211,221,611), \quad \varkappa_1=(J_0,K_1,K_3), \\
\tau_2&=(22,22,31,31,61,61;521), \quad \varkappa_2=(H_1,H_1,H_2,H_2,H_3,H_3;J_0), \quad \tau_4=(51).
\end{aligned}
\end{equation}
Then the \(2\)-class field tower of \(K\) has length \(\ell_2(K)=2\)
and metabelian Galois group \(\mathrm{Gal}(\mathrm{F}_2^\infty(K)/K)\)
isomorphic to the periodic group \(\langle 512,1605\rangle-\#1;2\) (see \S\
\ref{ss:FirstTree})
with order \(1024\) and coclass \(4\).
\end{theorem}


\begin{example}
\label{exm:ExcitedState}
Only a single imaginary quadratic fields \(K\) with \(2\)-class group of type \((4,4)\)
and Artin pattern in Formula
\eqref{eqn:APES}
is known up to now.
The discriminant \(d_K\) of this field is
\(-14\,547=(-3)\cdot 13\cdot 373\).
\end{example}


\begin{proof}
(Proof of Theorem
\ref{thm:ExcitedState})
None of the \(2\)-groups in the SmallGroups library
\cite{BEO}
hits the pattern up to the first layer \(\tau_0=(22)\), \(\tau_1=(211,221,611)\),
let alone the complete pattern up to the second layer \(\tau_2=(22,22,31,31,61,61;521)\).
However, the capitulation types \(\varkappa_1\) and \(\varkappa_2\) in Formula
\eqref{eqn:APES}
and Formula
\eqref{eqn:APGS}
coincide precisely,
whereas the abelian quotient invariants in these formulas
coincide only partially, while some components show an increase of
\(\tau_1=(211,221,411)\) and \(\tau_2=(22,22,31,31,41,41;321)\).
As explained in section
\ref{ss:FirstTree},
these observations suggest that the Artin pattern
\eqref{eqn:APES}
is an excited state of the ground state
\eqref{eqn:APGS}.
Hence, \(\langle 256,284\rangle\) is selected as starting point for extending the pool
by \(2\)-groups with orders bigger than \(512\), and thus outside of the SmallGroups database.
Descendants are generated with depth two (including grand children)
but only step size \(1\) (since bifurcations to coclass higher than \(4\) cannot occur).
Among the \(18\) groups generated in this manner,
there are \(4\) hits \(\langle 512,1605\rangle-\#1;i\) with \(1\le i\le 4\)
of the complete pattern \(\tau_0\), \(\tau_1\) and \(\tau_2\).
But only \(\langle 512,1605\rangle-\#1;2\) has correct capitulation type
(Table
\ref{tbl:Variants}).
\end{proof}


\begin{theorem}
\label{thm:OtherTreeI} (Another tree.)
Let \(K\) be an imaginary quadratic field
with Artin pattern
\begin{equation}
\label{eqn:APOtherTreeI}
\begin{aligned}
\tau_0&=(22), \quad \tau_1=(211,311,511), \quad \varkappa_1=(J_0,K_2,K_2), \\
\tau_2&=(31,31,31,31,51,51;421), \quad \varkappa_2=(H_1,H_1,H_2,H_2,H_2,H_2;H_2), \quad \tau_4=(41).
\end{aligned}
\end{equation}
Then the \(2\)-class field tower of \(K\) has length \(\ell_2(K)=2\)
and metabelian Galois group \(\mathrm{Gal}(\mathrm{F}_2^\infty(K)/K)\)
isomorphic to the periodic group \(\langle 512,1598\rangle\) with centre \(\zeta\simeq C_4\times C_2\)
on the tree \(\mathcal{T}^4(\langle 128,33\rangle)\).
\end{theorem}


\begin{example}
\label{exm:OtherTreeI}
Only a single imaginary quadratic fields \(K\) with \(2\)-class group of type \((4,4)\)
and Artin pattern in Formula
\eqref{eqn:APOtherTreeI}
is known up to now.
The discriminant \(d_K\) of this field is
\(-17\,476=(-4)\cdot 17\cdot 257\).
\end{example}


\renewcommand{\arraystretch}{1.1}

\begin{table}[ht]
\caption{Comparison of the candidates for Artin pattern \eqref{eqn:APOtherTreeI}}
\label{tbl:Candidates}
\begin{center}
{\small
\begin{tabular}{|c||c|c|c|}
\hline
 Vertex                      & \(\varkappa_1\) & \(\varkappa_2\)                    & \(\tau_4\) \\
\hline
 \(\langle 512,1595\rangle\) & \(J_0,K_2,J_0\) & \(H_1,H_1,H_2,H_2,G,G;H_2\)        & \((41)\)   \\
 \(\langle 512,1596\rangle\) & \(J_0,K_2,K_3\) & \(H_1,H_1,H_2,H_2,H_3,H_3;J_0\)    & \((41)\)   \\
 \(\langle 512,1597\rangle\) & \(J_0,K_2,K_1\) & \(H_1,H_1,H_2,H_2,H_1,H_1;J_0\)    & \((41)\)   \\
 \(\langle 512,1598\rangle\) & \(J_0,K_2,K_2\) & \(H_1,H_1,H_2,H_2,H_2,H_2;H_2\)    & \((41)\)   \\
\hline
\end{tabular}
}
\end{center}
\end{table}


\begin{proof}
(Proof of Theorem
\ref{thm:OtherTreeI})
Searching the SmallGroups database
\cite{BEO}
for the pattern \(\tau_0=(22)\), \(\tau_1=(211,311,511)\) and \(\tau_2=(31,31,31,31,51,51;421)\)
yields \(4\) groups \(\langle 512,i\rangle\) with identifiers \(1595\le i\le 1598\).
But only \(\langle 512,1598\rangle\) possesses the correct capitulation type
\(\varkappa_1=(J_0,K_3,K_3)\) with two equal kernels
and \(\varkappa_2=(H_2,H_2,H_3,H_3,H_2,H_2;H_2)\).
See Table
\ref{tbl:Candidates}.
Here, we have \(\mathrm{Cl}_2(N_1)=(211)\), \(\mathrm{Cl}_2(N_{1k})=(31)\), for \(k\in\lbrace 1,2\rbrace\),
and indeed \(\tau_4=(41)\)
\cite[Cor. 3.2]{BeSn}.
\end{proof}


\subsection{The first descendant tree of \(2\)-groups with type \((4,4)\)}
\label{ss:FirstTree}
\noindent
The first descendant tree \(\mathcal{T}^4(R)\) of finite \(2\)-groups \(G\)
with abelian type invariants \((4,4)\) and fixed coclass \(\mathrm{cc}(G)=4\)
was discovered in May 2019 by analyzing the second \(2\)-class group
\(\mathrm{Gal}(\mathrm{F}_2^2(K)/K)\) of imaginary quadratic fields \(K\)
with discriminants \(d_K\in\lbrace -8\,103,\ -11\,972,\ -14\,547\rbrace\)
and \(2\)-class group \(\mathrm{Cl}_2(K)\simeq C_4\times C_4\).
Its metabelian root \(R=\langle 128,35\rangle\) is connected with the
general abelian root \(A=\langle 16,2\rangle\simeq C_4\times C_4\) by the path
\(R\stackrel{2}{\to}\langle 32,2\rangle\to A\).
Up to order \(1024\), the tree is drawn in Figure
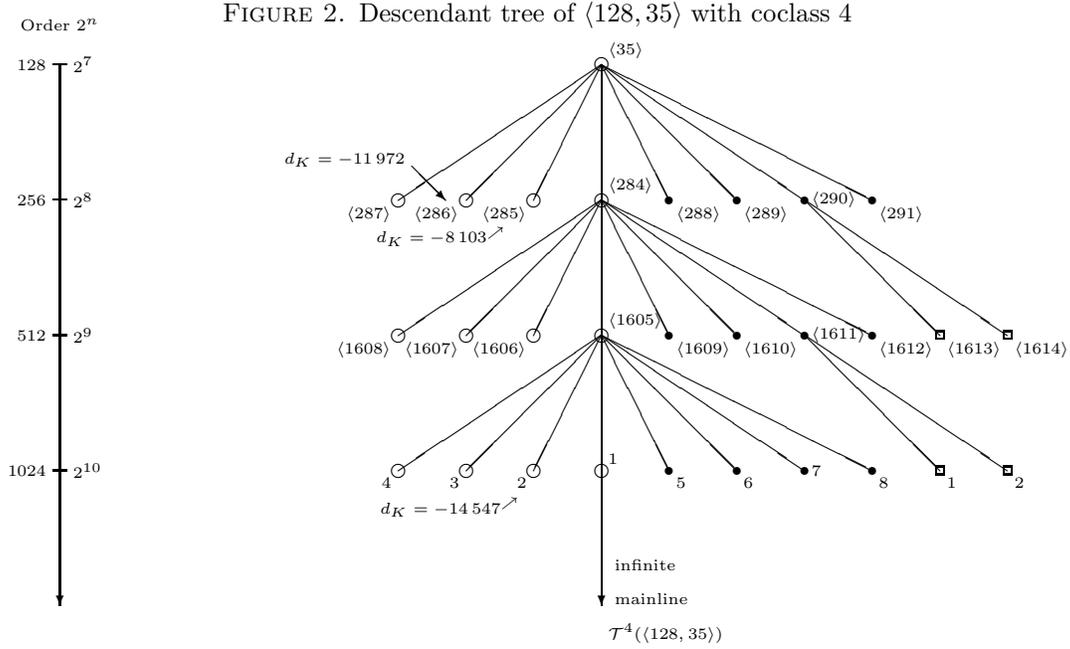
\ref{fig:FirstTree}.


\begin{figure}[ht]
\caption{Descendant tree of \(\langle 128,35\rangle\) with coclass \(4\)}
\label{fig:FirstTree}

{\tiny

\setlength{\unitlength}{0.9cm}
\begin{picture}(10,9)(-6,-8.5)

\put(-8,0.5){\makebox(0,0)[cb]{Order \(2^n\)}}
\put(-8,0){\line(0,-1){6}}
\multiput(-8.1,0)(0,-2){4}{\line(1,0){0.2}}
\put(-8.2,0){\makebox(0,0)[rc]{\(128\)}}
\put(-7.8,0){\makebox(0,0)[lc]{\(2^7\)}}
\put(-8.2,-2){\makebox(0,0)[rc]{\(256\)}}
\put(-7.8,-2){\makebox(0,0)[lc]{\(2^8\)}}
\put(-8.2,-4){\makebox(0,0)[rc]{\(512\)}}
\put(-7.8,-4){\makebox(0,0)[lc]{\(2^9\)}}
\put(-8.2,-6){\makebox(0,0)[rc]{\(1024\)}}
\put(-7.8,-6){\makebox(0,0)[lc]{\(2^{10}\)}}
\put(-8,-6){\vector(0,-1){2}}

\multiput(0,0)(0,-2){4}{\circle{0.2}}
\multiput(0,0)(0,-2){3}{\line(0,-1){2}}
\multiput(-3,-2)(0,-2){3}{\circle{0.2}}
\multiput(-2,-2)(0,-2){3}{\circle{0.2}}
\multiput(-1,-2)(0,-2){3}{\circle{0.2}}
\multiput(1,-2)(0,-2){3}{\circle*{0.1}}
\multiput(2,-2)(0,-2){3}{\circle*{0.1}}
\multiput(3,-2)(0,-2){3}{\circle*{0.1}}
\multiput(4,-2)(0,-2){3}{\circle*{0.1}}
\multiput(0,0)(0,-2){3}{\line(-3,-2){3}}
\multiput(0,0)(0,-2){3}{\line(-1,-1){2}}
\multiput(0,0)(0,-2){3}{\line(-1,-2){1}}
\multiput(0,0)(0,-2){3}{\line(1,-2){1}}
\multiput(0,0)(0,-2){3}{\line(1,-1){2}}
\multiput(0,0)(0,-2){3}{\line(3,-2){3}}
\multiput(0,0)(0,-2){3}{\line(2,-1){4}}
\multiput(4.95,-4.05)(0,-2){2}{\framebox(0.1,0.1){}}
\multiput(5.95,-4.05)(0,-2){2}{\framebox(0.1,0.1){}}
\multiput(3,-2)(0,-2){2}{\line(1,-1){2}}
\multiput(3,-2)(0,-2){2}{\line(3,-2){3}}

\put(0.1,0.1){\makebox(0,0)[lb]{\(\langle 35\rangle\)}}

\put(-2.9,-1.3){\makebox(0,0)[rt]{\(d_K=-11\,972\)}}
\put(-2.8,-1.5){\vector(1,-1){0.5}}
\put(-3.1,-2.1){\makebox(0,0)[rt]{\(\langle 287\rangle\)}}
\put(-2.1,-2.1){\makebox(0,0)[rt]{\(\langle 286\rangle\)}}
\put(-1.1,-2.1){\makebox(0,0)[rt]{\(\langle 285\rangle\)}}
\put(-1.4,-2.4){\makebox(0,0)[rt]{\(d_K=-8\,103^\nearrow\)}}

\put(0.1,-1.9){\makebox(0,0)[lb]{\(\langle 284\rangle\)}}

\put(1.1,-2.1){\makebox(0,0)[lt]{\(\langle 288\rangle\)}}
\put(2.1,-2.1){\makebox(0,0)[lt]{\(\langle 289\rangle\)}}
\put(3.1,-2.0){\makebox(0,0)[lc]{\(\langle 290\rangle\)}}
\put(4.1,-2.1){\makebox(0,0)[lt]{\(\langle 291\rangle\)}}

\put(-3.1,-4.1){\makebox(0,0)[rt]{\(\langle 1608\rangle\)}}
\put(-2.1,-4.1){\makebox(0,0)[rt]{\(\langle 1607\rangle\)}}
\put(-1.1,-4.1){\makebox(0,0)[rt]{\(\langle 1606\rangle\)}}

\put(0.1,-3.9){\makebox(0,0)[lb]{\(\langle 1605\rangle\)}}

\put(1.1,-4.1){\makebox(0,0)[lt]{\(\langle 1609\rangle\)}}
\put(2.1,-4.1){\makebox(0,0)[lt]{\(\langle 1610\rangle\)}}
\put(3.1,-4.0){\makebox(0,0)[lc]{\(\langle 1611\rangle\)}}
\put(4.1,-4.1){\makebox(0,0)[lt]{\(\langle 1612\rangle\)}}
\put(5.1,-4.1){\makebox(0,0)[lt]{\(\langle 1613\rangle\)}}
\put(6.1,-4.1){\makebox(0,0)[lt]{\(\langle 1614\rangle\)}}

\put(-3.1,-6.1){\makebox(0,0)[rt]{\(4\)}}
\put(-2.1,-6.1){\makebox(0,0)[rt]{\(3\)}}
\put(-1.1,-6.1){\makebox(0,0)[rt]{\(2\)}}
\put(-1.2,-6.4){\makebox(0,0)[rt]{\(d_K=-14\,547^\nearrow\)}}

\put(0.1,-5.9){\makebox(0,0)[lb]{\(1\)}}

\put(1.1,-6.1){\makebox(0,0)[lt]{\(5\)}}
\put(2.1,-6.1){\makebox(0,0)[lt]{\(6\)}}
\put(3.1,-6.0){\makebox(0,0)[lc]{\(7\)}}
\put(4.1,-6.1){\makebox(0,0)[lt]{\(8\)}}
\put(5.1,-6.1){\makebox(0,0)[lt]{\(1\)}}
\put(6.1,-6.1){\makebox(0,0)[lt]{\(2\)}}

\put(0,-6){\vector(0,-1){2}}
\put(0.2,-7.4){\makebox(0,0)[lc]{infinite}}
\put(0.2,-7.9){\makebox(0,0)[lc]{mainline}}
\put(1.8,-8.4){\makebox(0,0)[rc]{\(\mathcal{T}^4(\langle 128,35\rangle)\)}}

\end{picture}
}
\end{figure}


The mainline is characterized graph theoretically by the descendant numbers \((N,C)=(8,2)\).
It is possible to partition the eight immediate descendants of a mainline vertex
group theoretically:
Four vertices (\(\bigcirc\)) possess a centre \(\zeta\simeq C_4\times C_2\),
one of them belongs to the mainline, the other three are terminal.
The other four vertices \((\bullet\)) have \(\zeta\simeq C_2\times C_2\),
the third of them has nuclear rank \(\nu=1\) and descendant numbers \((N,C)=(2,0)\),
the other three are leaves.
Both vertices of depth two (\(\square\)) have \(\zeta\simeq C_2\).
(By \(\zeta=\zeta_1\) we denote the first term of the upper central series.)

The tree is an excellent paradigm for the principle of
\textit{stabilization} and \textit{polarization} of certain components of the Artin pattern.
We avoid complications by restricting ourselves
to the vertices (\(\bigcirc\)) with centre \(\zeta\simeq C_4\times C_2\)
which are located on and left from the mainline in Figure
\ref{fig:FirstTree}.

Among the transfer \textit{targets} (abelian quotient invariants),
the first two components of the first layer \(\tau_1=(211,221,\ast)\) and
the first four components of the second layer \(\tau_2=(22,22,31,31,\ast,\ast;\ast)\)
remain stable on the entire tree.
Similarly for the transfer \textit{kernels} (capitulation kernels),
the first two components of the first layer \(\varkappa_1=(J_0,K_1,\ast)\) and
the first four components of the second layer \(\varkappa_2=(H_1,H_1,H_2,H_2,\ast,\ast;\ast)\)
remain stable.

The polarization appears in a different manner for targets and kernels.
The polarized components of the transfer \textit{targets}
become strictly bigger from parent to descendant (according to their isotonic monotony principle),
i.e., with increasing order.
Starting at the root, we successively have 
\(\tau_1(3)\in\lbrace 311,411,511,611,\ldots\rbrace\),
\(\tau_2(5)=\tau_2(6)\in\lbrace 31,41,51,61,\ldots\rbrace\),
\(\tau_2(7)\in\lbrace 221,321,421,521,\ldots\rbrace\), and
\(\tau_4\in\lbrace 21,31,41,51,\ldots\rbrace\).
This phenomenon gives rise to countably many \textit{states},
corresponding to the infinitely many branches of the tree in vertical direction
(which are all isomorphic as graphs for this tree with period length \(1\)),
and indicated in Table
\ref{tbl:States}.


\renewcommand{\arraystretch}{1.1}

\begin{table}[ht]
\caption{States with increasing polarized transfer targets}
\label{tbl:States}
\begin{center}
\begin{tabular}{|c||r|c||c|c|c|c|}
\hline
 State       & Order    & Identifiers             & \(\tau_1(3)\) & \(\tau_2(5)\) & \(\tau_2(7)\) & \(\tau_4\) \\
\hline
 Root        &  \(128\) & \(35\)                  & \((311)\)     & \((31)\)      & \((221)\)     & \((21)\)   \\
\hline
 Ground      &  \(256\) & \(284,285,286,287\)     & \((411)\)     & \((41)\)      & \((321)\)     & \((31)\)   \\
 1st Excited &  \(512\) & \(1605,1606,1607,1608\) & \((511)\)     & \((51)\)      & \((421)\)     & \((41)\)   \\
 2nd Excited & \(1024\) & \(1,2,3,4\)             & \((611)\)     & \((61)\)      & \((521)\)     & \((51)\)   \\
\hline
\end{tabular}
\end{center}
\end{table}


On the other hand,
the polarized components of the transfer \textit{kernels}
(which may become smaller from parent to descendant, according to their antitonic monotony principle)
vary in the same way for each fixed order
\(\varkappa_1(3)\in\lbrace J_0,K_3,K_1,K_2\rbrace\),
\(\varkappa_2(5)=\varkappa_2(6)\in\lbrace G,H_3,H_1,H_2\rbrace\), and
\(\varkappa_2(7)\in\lbrace H_1,J_0,H_1,J_0\rbrace\),
with respect to increasing SmallGroups identifiers.
Hence, there are only finitely many \textit{variants} within a fixed state in horizontal direction,
characterized by distinct capitulation types in Table
\ref{tbl:Variants}.


\renewcommand{\arraystretch}{1.1}

\begin{table}[ht]
\caption{Variants with distinct polarized transfer kernels}
\label{tbl:Variants}
\begin{center}
\begin{tabular}{|c||c||c|c|c|}
\hline
 Variant  & Identifiers                                                                & \(\varkappa_1(3)\) & \(\varkappa_2(5)\) & \(\varkappa_2(7)\) \\
\hline
 Mainline & \(\langle 128,35\rangle,\langle 256,284\rangle,\langle 512,1605\rangle,1\) & \(J_0\)            & \(G\)              & \(H_1\)            \\
\hline
 First    & \(\langle 256,285\rangle,\langle 512,1606\rangle,2\)                       & \(K_3\)            & \(H_3\)            & \(J_0\)            \\
 Second   & \(\langle 256,286\rangle,\langle 512,1607\rangle,3\)                       & \(K_1\)            & \(H_1\)            & \(H_1\)            \\
 Third    & \(\langle 256,287\rangle,\langle 512,1608\rangle,4\)                       & \(K_2\)            & \(H_2\)            & \(J_0\)            \\
\hline
\end{tabular}
\end{center}
\end{table}


The tree is also an excellent paradigm for the \textit{mainline principle},
which states that mainline vertices have bigger capitulation kernels in their polarization.
Indeed, we have
\(J_0>K_i\) for \(\varkappa_1(3)\),
\(G>H_i\) for \(\varkappa_2(5)=\varkappa_2(6)\), and 
\(H_1>J_0\) for \(\varkappa_2(7)\).


\subsection{Non-metabelian towers of Hilbert \(2\)-class fields}
\label{ss:Length3}
\noindent
Table
\ref{tbl:Imaginary3}
shows some discriminants \(0>d_K>-10^5\) of imaginary quadratic fields \(K\)
with \(\mathrm{Cl}_2(K)\simeq C_4\times C_4\),
their factorization into three prime discriminants
and the metabelianization \(G/G^{\prime\prime}\simeq\mathrm{Gal}(\mathrm{F}_2^2(K)/K)\)
of the Galois group \(G=\mathrm{Gal}(\mathrm{F}_2^\infty(K)/K)\)
of the non-metabelian \(2\)-class field tower of \(K\)
in the notation of the SmallGroups library
\cite{BEO}
and ANUPQ package
\cite{GNO}.


\renewcommand{\arraystretch}{1.1}

\begin{table}[ht]
\caption{Discriminants of imaginary quadratic fields \(K\) with \(\ell_2(K)=3\)}
\label{tbl:Imaginary3}
\begin{center}
\begin{tabular}{|r|c||c|c|}
\hline
 \(d_K\)      & Factorization             & \(\mathrm{Gal}(\mathrm{F}_2^2(K)/K)\) & \(\tau_4\) \\
\hline
  \(-6\,123\) & \((-3)\cdot 13\cdot 157\) & \(\langle 512,227\rangle-\#2;140\)    & \((421)\)  \\
 \(-25\,144\) & \(8\cdot (-7)\cdot 449\)  & \(\langle 512,231\rangle-\#1;4\)      & \((321)\)  \\
\hline
\end{tabular}
\end{center}
\end{table}

\noindent
The root path of the occurring metabelian \(2\)-groups with orders \(2048\) and \(1024\) is given by
\eqref{eqn:RootPathI3}.

\begin{equation}
\label{eqn:RootPathI3}
\begin{aligned}
\#2;140\stackrel{2}{\to}\langle 512,227\rangle\stackrel{2}{\to}\langle 128,26\rangle\stackrel{2}{\to}\langle 32,2\rangle\to\langle 16,2\rangle \\
\#1;4\to\langle 512,231\rangle\stackrel{2}{\to}\langle 128,26\rangle\stackrel{2}{\to}\langle 32,2\rangle\to\langle 16,2\rangle
\end{aligned}
\end{equation}

\noindent
In accordance with the Shafarevich criterion,
these metabelianizations have relation rank \(d_2(G/G^{\prime\prime})=4\),
and therefore cannot be \(2\)-tower groups of imaginary quadratic fields.
The coclass is \(6\), resp. \(5\).


\begin{theorem}
\label{thm:ThreeStagesI1} (Three-stage tower.)
Let \(K\) be an imaginary quadratic field
with Artin pattern
\begin{equation}
\label{eqn:APThreeStagesI1}
\begin{aligned}
\tau_0&=(22), \quad \tau_1=(211,221,411), \quad \varkappa_1=(J_0,K_3,K_1), \\
\tau_2&=(22,22,41,41,222,2111;3111), \quad \varkappa_2=(H_3,H_3,H_1,H_1,H_2,H_1;J_0), \quad \tau_4=(321).
\end{aligned}
\end{equation}
The \(2\)-class tower of \(K\) has length \(\ell_2(K)=3\)
and the metabelianization \(G/G^{\prime\prime}=\mathrm{Gal}(\mathrm{F}_2^2(K)/K)\)
of the Galois group \(G=\mathrm{Gal}(\mathrm{F}_2^\infty(K)/K)\)
is isomorphic to \(\langle 512,231\rangle-\#1;4\) with order \(1024\).
The group \(G\) itself is isomorphic to \(\langle 512,231\rangle-\#2;j\) with some \(j\in\lbrace 36,39,52,55\rbrace\)
and order \(2048\).
\end{theorem}


\begin{example}
\label{exm:ThreeStagesI1}
There are only two known discriminants \(0>d_K>-10^5\)
of imaginary quadratic fields \(K\) with \(\mathrm{Cl}_2(K)\simeq C_4\times C_4\)
and Artin pattern in Formula
\eqref{eqn:APThreeStagesI1}:
\begin{equation}
\label{eqn:ExmThreeStagesI1}
\begin{aligned}
-25\,144&=8\cdot (-7)\cdot 449, \qquad
&-37\,407&=(-3)\cdot 37\cdot 337.
\end{aligned}
\end{equation}
\end{example}


\begin{proof}
(Proof of Theorem
\ref{thm:ThreeStagesI1})
Searching the SmallGroups database
\cite{BEO}
for the pattern \(\tau_0=(22)\) and \(\tau_1=(211,221,411)\), up to the first layer alone,
yields \(8\) groups \(\langle 256,i\rangle\) with \(246\le i\le 249\) and \(284\le i\le 287\)
and \(28\) groups \(\langle 512,j\rangle\) with \(227\le j\le 242\), \(1539\le j\le 1546\) and \(1609\le j\le 1612\).
But none of the groups in the SmallGroups library hits the complete pattern up to the second layer,
including \(\tau_2=(22,22,41,41,222,2111;3111)\).
Some of the \(36\) candidates have too big components in their second layer pattern \(\tau_2\):
The four groups \(\langle 256,i\rangle\) with \(284\le i\le 287\) have inadequate Frattini invariants \((321)\) instead of \((3111)\),
the eight groups \(\langle 512,j\rangle\) with \(1539\le j\le 1542\) and \(1609\le j\le 1612\) have two inadequate invariants \((51)\) instead of \((41)\),
and they can be eliminated as parents for extending the pool by the \(p\)-group generation algorithm
\cite{Nm,Ob},
according to the \textit{monotony principle}
\cite[\S\ 1.7, Thm. 1.21, p. 79]{Ma2016b}.
The four groups \(\langle 512,j\rangle\) with \(1543\le j\le 1546\)
can be eliminated, because they are terminal (without descendants).

Since all the other groups share \(\langle 128,26\rangle\) as their common parent,
this group is chosen as starting point for selectively constructing an extended pool
of \(2\)-groups with orders bigger than \(512\), and thus outside of the SmallGroups database.
Descendants are generated with depth two (including grand children)
and step sizes \(1\) and \(2\) (including bifurcations to coclass higher than \(4\)).
Among the \(1754\) groups generated in this manner,
there are \(176\) hits of the complete pattern \(\tau_0=(22)\), \(\tau_1=(211,221,411)\),
\(\tau_2=(22,22,41,41,222,2111;3111)\) and \(\tau_4=(321)\) of abelian type invariants,
but without considering any capitulation patterns \(\varkappa_1\) and \(\varkappa_2\).
\(16\) of the hits are metabelian with order \(1024\) and coclass \(5\),
and can be partitioned into four batches with \(4\) members each.
The remaining \(160\) hits are non-metabelian with derived length \(3\), order \(2048\) and coclass \(6\),
and can be partitioned into four batches, one with \(64\) members and three with \(32\) members each.
There is a one-to-one correspondence between metabelian and non-metabelian batches,
where the members of the former are metabelianizations of the latter,
and furthermore both share a common parent.
The pair of batches with parent \(\langle 512,227\rangle\) is eliminated,
due to an inadequate capitulation type \(\varkappa_1=(J_0,K_3,J_0)\).
(For this pair of batches, \(16\) non-metabelian groups map to the same metabelianization,
for the other pairs of batches, only \(8\) map to the same.)
Similarly, the pair of batches with parent \(\langle 512,228\rangle\) is eliminated,
due to an inadequate capitulation type \(\varkappa_1=(J_0,K_3,K_3)\).
The two pairs of batches with parents \(\langle 512,231\rangle\), resp. \(\langle 512,232\rangle\),
have admissible \(\varkappa_1=(J_0,K_3,K_1)\), resp. \(\varkappa_1=(J_0,K_3,K_2)\).
For the final decision, the second layer capitulation type \(\varkappa_2\) must be taken into acount.
The pair of batches with parent \(\langle 512,232\rangle\) is eliminated,
because \(\varkappa_2=(H_3,H_3,H_2,H_2,\ast,H_1;J_0)\) assigns different capitulation kernels
to the components with abelian type invariants \((41)\) and \((2111)\).
For the last remaining pair of batches with parent \(\langle 512,231\rangle\), Table
\ref{tbl:PolThreeStagesI1}
shows the development of the Artin pattern along the root path,
for each of the four members \(\#1;i\) with \(i\in\lbrace 3,4,7,8\rbrace\) of the metabelian batch.
They all share the common \textit{stabilization}
\(\varkappa_2=(H_3,H_3,H_1,H_1,\ast,H_1;J_0)\),
but only \(\#1;4\) possesses the correct \textit{polarization} with \(\ast\) replaced by \(H_2\)
such that the component with abelian type invariants \((222)\)
has a capitulation kernel different from all other kernels.


\renewcommand{\arraystretch}{1.0}

\begin{table}[ht]
\caption{Development of the Artin pattern along the root path}
\label{tbl:PolThreeStagesI1}
\begin{center}
{\small
\begin{tabular}{|c||c|c|c|c|c|}
\hline
 Vertex                     & \(\tau_1\)      & \(\varkappa_1\) & \(\tau_2\)                    & \(\varkappa_2\)                 & \(\tau_4\) \\
\hline
 \(\langle 16,2\rangle\)    & \(21,21,21\)    & \(J_0,J_0,J_0\) & \(2,2,2,2,2,2;11\)            & \(G,G,G,G,G,G;G\)               & \((0)\)    \\
 \(\langle 32,2\rangle\)    & \(211,211,211\) & \(J_0,J_0,J_0\) & \(21,21,21,21,21,21;111\)     & \(G,G,G,G,G,G;G\)               & \((1)\)    \\
 \(\langle 128,26\rangle\)  & \(211,221,311\) & \(J_0,K_3,J_0\) & \(22,22,31,31,211,211;2111\)  & \(H_3,H_3,G,G,G,G;H_3\)         & \((21)\)   \\
 \(\langle 512,231\rangle\) & \(211,221,411\) & \(J_0,K_3,K_1\) & \(22,22,41,41,221,2111;3111\) & \(H_3,H_3,H_1,H_1,G,H_1;J_0\)   & \((311)\)  \\
 \(\#1;3\)                  & \(211,221,411\) & \(J_0,K_3,K_1\) & \(22,22,41,41,222,2111;3111\) & \(H_3,H_3,H_1,H_1,H_1,H_1;J_0\) & \((321)\)  \\
 \(\#1;4\)                  & \(211,221,411\) & \(J_0,K_3,K_1\) & \(22,22,41,41,222,2111;3111\) & \(H_3,H_3,H_1,H_1,H_2,H_1;J_0\) & \((321)\)  \\
 \(\#1;7\)                  & \(211,221,411\) & \(J_0,K_3,K_1\) & \(22,22,41,41,222,2111;3111\) & \(H_3,H_3,H_1,H_1,G,H_1;J_0\)   & \((321)\)  \\
 \(\#1;8\)                  & \(211,221,411\) & \(J_0,K_3,K_1\) & \(22,22,41,41,222,2111;3111\) & \(H_3,H_3,H_1,H_1,H_3,H_1;J_0\) & \((321)\)  \\
\hline
\end{tabular}
}
\end{center}
\end{table}


Among the eight members of the associated non-metabelian batch,
four are capable with relation rank \(d_2=4\),
and only the other four siblings \(\#2;j\) with \(j\in\lbrace 36,39,52,55\rbrace\)
which are terminal with required relation rank \(d_2=3\)
are candidates for the \(2\)-tower group \(G\) with three stages.

Here, we have \(\mathrm{Cl}_2(N_1)=(211)\), but
\(\mathrm{Cl}_2(N_{11})=(222)\), \(\mathrm{Cl}_2(N_{12})=(2111)\) with order bigger than \(16\), 
and thus \(\tau_4=(321)\) with rank three
\cite[Cor. 3.2]{BeSn}.
\end{proof}


\begin{theorem}
\label{thm:ThreeStagesI2} (Tower with \(\ell_2\ge 3\).)
Let \(K\) be an imaginary quadratic field
with Artin pattern
\begin{equation}
\label{eqn:APThreeStagesI2}
\begin{aligned}
\tau_0&=(22), \quad \tau_1=(211,221,511), \quad \varkappa_1=(J_0,K_1,K_3), \\
\tau_2&=(22,22,51,51,222,2111;4111), \quad \varkappa_2=(H_1,H_1,H_2,H_2,H_3,H_2;J_0), \quad \tau_4=(421).
\end{aligned}
\end{equation}
The \(2\)-class tower of \(K\) has length \(\ell_2(K)\ge 3\)
and the metabelianization \(G/G^{\prime\prime}=\mathrm{Gal}(\mathrm{F}_2^2(K)/K)\)
of the Galois group \(G=\mathrm{Gal}(\mathrm{F}_2^\infty(K)/K)\)
is isomorphic to \(\langle 512,227\rangle-\#2;140\) of order \(2048\) and coclass \(6\).
The smallest possibility for the \(2\)-tower group \(G\) itself
is an isomorphism to one of the four non-metabelian groups
\(\langle 512,227\rangle-\#3;j\) with \(j\in\lbrace 39,62,231,254\rbrace\) (see \S\
\ref{ss:LimitsAndIPAD})
of order \(4096\), coclass \(7\) and derived length \(3\).
\end{theorem}


\begin{example}
\label{exm:ThreeStagesI2}
\(-6\,123\) is the single known discriminant \(0>d_K>-10^5\)
of an imaginary quadratic field \(K\) with \(\mathrm{Cl}_2(K)\simeq C_4\times C_4\)
and Artin pattern in Formula
\eqref{eqn:APThreeStagesI2}.
It is mentioned in
\cite[Exm. 7.2, p. 1192]{BeSn}
without explicit invariant \(\tau_4=(421)\).
However, \(\ell_2(K)\) and \(G\) are unknown (\S\
\ref{ss:LimitsAndIPAD}).
\end{example}


\begin{proof}
(Proof of Theorem
\ref{thm:ThreeStagesI2})
We have intentionally delayed this proof,
because we shall need our experience from the proof of Theorem
\ref{thm:ThreeStagesI1}.
Searching the SmallGroups database
\cite{BEO}
for the pattern \(\tau_0=(22)\) and \(\tau_1=(211,221,511)\), up to the first layer alone,
yields \(8\) groups \(\langle 512,j\rangle\) with
\(1535\le j\le 1538\) and parent \(\langle 256,246\rangle\), resp.
\(1605\le j\le 1608\) and parent \(\langle 256,284\rangle\).
But none of the groups in the SmallGroups library hits the complete pattern up to the second layer,
including \(\tau_2=(22,22,51,51,222,2111;4111)\).
Construction of descendants with starting point \(\langle 256,246\rangle\) or \(\langle 256,284\rangle\)
can be terminated after depth one already,
according to the \textit{monotony principle}
\cite[\S\ 1.7, Thm. 1.21, p. 79]{Ma2016b},
since certain components of \(\tau_2\) become too big.
At this point we lack any further traces.

Guided by the proof of Theorem
\ref{thm:ThreeStagesI1},
we remember the promising starting group \(\langle 128,26\rangle\)
with \(\tau_1=(211,221,311)\) and \(\tau_2=(22,22,31,31,211,211;2111)\),
which would be compatible with our search pattern
\eqref{eqn:APThreeStagesI1},
according to the monotony principle.
Indeed, a similar construction of \(1754\) descendants as in the proof of Theorem
\ref{thm:ThreeStagesI1}
yields \(16\) hits by metabelian groups of order \(2048\)
with common parent \(\langle 512,227\rangle\).
Only one of them, namely \(\langle 512,227\rangle-\#2;140\) possesses the required capitulation type with
\(\varkappa_1=(J_0,K_1,K_3)\) and \(\varkappa_2=(H_1,H_1,H_2,H_2,H_3,H_2;J_0)\).

Finally, we construct the \(1016\) immediate descendants of all step sizes \(1,2,3\)
of \(\langle 512,227\rangle\) (including the trifurcation to coclass \(7\))
and find eight non-metabelian groups with metabelianization \(\langle 512,227\rangle-\#2;140\).
Four of them are capable with relation rank \(d_2=4\),
and only the other four siblings \(\#3;j\) with \(j\in\lbrace 39,62,231,254\rbrace\)
which are terminal with required relation rank \(d_2=3\)
are candidates for the \(2\)-tower group \(G\) with three stages.

Here, we have \(\mathrm{Cl}_2(N_1)=(211)\), but
\(\mathrm{Cl}_2(N_{11})=(222)\), \(\mathrm{Cl}_2(N_{12})=(2111)\) with order bigger than \(16\), 
and thus \(\tau_4=(421)\) with rank three
\cite[Cor. 3.2]{BeSn}.
\end{proof}


The situations in Theorem
\ref{thm:ThreeStagesI1}
resp. Theorem
\ref{thm:ThreeStagesI2}
were called \textit{advanced tree topologies} of type \textit{fork and siblings} in
\cite[\S\ 5, pp. 89--92]{Ma2016b}.
The fork is
\(\langle 512,231\rangle\) resp. \(\langle 512,227\rangle\).


\subsection{Limits and higher abelian quotient invariants}
\label{ss:LimitsAndIPAD}
\noindent
Imaginary quadratic fields \(K\) with Artin pattern in Formula
\eqref{eqn:APThreeStagesI2},
in particular the field with discriminant \(d_K=-6\,123\) in Example
\ref{exm:ThreeStagesI2},
give rise to problems with determining the exact length \(\ell_2(K)\ge 3\)
of their \(2\)-class field tower.
The following theorem illuminates the minimal candidates
for the \(2\)-tower group \(G=\mathrm{Gal}(\mathrm{F}_2^\infty(K)/K)\)
of such fields \(K\) in Theorem
\ref{thm:ThreeStagesI2}
more closely
from the viewpoint of finitely presented \(2\)-groups,
and shows that the tower length \(\ell_2(K)\) cannot be bounded by \(3\).


\begin{theorem}
\label{thm:Limits}
(Pro-\(2\) groups as limit groups) \\
The eight non-metabelian \(2\)-groups of order \(4096\)
with metabelianization \(M:=\langle 512,227\rangle-\#2;140\),
\begin{equation}
\label{eqn:Metabelian}
M=\langle x,y\mid x^4=\lbrack t_3,y,y\rbrack,\ y^4=\lbrack t_3,y,xy\rbrack,\ \lbrack x^2,y\rbrack=1, w=1\rangle,
\end{equation}
are the class-\(5\) quotients
\begin{equation}
\label{eqn:Quotients}
\begin{aligned}
G_{000}/\gamma_6G_{000}&=\langle 512,227\rangle-\#3;38, \text{ infinitely capable},\ d_2=4, \\
G_{001}/\gamma_6G_{001}&=\langle 512,227\rangle-\#3;39, \text{ terminal},\ d_2=3, \\
G_{010}/\gamma_6G_{010}&=\langle 512,227\rangle-\#3;62, \text{ terminal},\ d_2=3, \\
G_{011}/\gamma_6G_{011}&=\langle 512,227\rangle-\#3;63, \text{ finitely capable},\ d_2=4, \\
G_{100}/\gamma_6G_{100}&=\langle 512,227\rangle-\#3;230, \text{ finitely capable},\ d_2=4, \\
G_{101}/\gamma_6G_{101}&=\langle 512,227\rangle-\#3;231, \text{ terminal},\ d_2=3, \\
G_{110}/\gamma_6G_{110}&=\langle 512,227\rangle-\#3;254, \text{ terminal},\ d_2=3, \\
G_{111}/\gamma_6G_{111}&=\langle 512,227\rangle-\#3;255, \text{ finitely capable},\ d_2=4,
\end{aligned}
\end{equation}
of the following finitely presented groups:
\begin{equation}
\label{eqn:Limits}
G_{efg}:=\langle x,y\mid \lbrack x^2,y\rbrack w^e=1,\ x^4w^f=\lbrack t_3,y,y\rbrack,\ y^4w^g=\lbrack t_3,y^2\rbrack\rangle
\text{ with } 0\le e,f,g\le 1, \\
\end{equation}
where \(s_2:=\lbrack y,x\rbrack\) is the main commutator, \(t_3:=\lbrack s_2,y\rbrack\),
and \(w\) denotes the group word
\(w(x,y):=\lbrack t_3,s_2\rbrack=\lbrack \lbrack y,x,y\rbrack,\lbrack y,x\rbrack\rbrack\).

There exists a non-metabelian \(2\)-group \(\mathfrak{G}_{20}\) of order \(2^{20}\) with derived length \(4\), relation rank \(3\) and coclass \(9\)
with metabelianization \(M\) and the following root path:

{\scriptsize
\begin{equation}
\label{eqn:RootPath}
\#1;1\to\#1;1\to\#1;1\to\#2;2\stackrel{2}{\to}\#2;1\stackrel{2}{\to}\#1;4\to\#3;230
\stackrel{3}{\to}\langle 512,227\rangle\stackrel{2}{\to}\langle 128,26\rangle\stackrel{2}{\to}\langle 32,2\rangle\to\langle 16,2\rangle.
\end{equation}}

\noindent
This group \(\mathfrak{G}_{20}\) is descendant \(G_{100}/\gamma_{12}G_{100}\) of \(G_{100}/\gamma_6G_{100}\)
and class-\(11\) quotient of \(G_{100}\).
\end{theorem}


\noindent
However, in spite of Theorem
\ref{thm:Limits},
we cannot determine the \(2\)-tower group \(G\) of the field \(K\) with discriminant \(d_K=-6\,123\),
because its abelian quotient invariants \(\tau_{2,1}\) of higher order in Table
\ref{tbl:IPAD2}
do not match any known finite \(2\)-group.


\renewcommand{\arraystretch}{1.1}

\begin{table}[ht]
\caption{Selected higher abelian quotient invariants of known groups and \(K\)}
\label{tbl:IPAD2}
\begin{center}
\begin{tabular}{|c|c|}
\hline
 Group resp. Field & \(\tau_{2,1}\) \\
\hline
 \(M=\langle 512,227\rangle-\#2;140\) & \(\lbrack(22;211,211,411)^2,(51;51,51,411)^2,(2111;211,\ldots),\ldots\rbrack\)  \\
 \(\langle 512,227\rangle-\#3;39\)    & \(\lbrack(22;221,221,411)^2,(51;52,52,411)^2,(2111;2111,\ldots),\ldots\rbrack\)  \\
 \(\mathfrak{G}_{20}\)                & \(\lbrack(22;221,221,411)^2,(51;52,52,411)^2,(2111;2211,\ldots),\ldots\rbrack\)  \\
\hline
 \(K=\mathbb{Q}(\sqrt{-6123})\)       & \(\lbrack(22;311,311,411)^2,(51;62,62,411)^2,(2111;3111,\ldots),\ldots\rbrack\)  \\
\hline
\end{tabular}
\end{center}
\end{table}

\subsection{Unexpected two-stage towers of Hilbert \(2\)-class fields}
\label{ss:UnexpLength2I}
\noindent
Table
\ref{tbl:UnexpImaginary2}
shows a discriminant \(0>d_K>-10^5\) of an imaginary quadratic field \(K\)
with \(\mathrm{Cl}_2(K)\simeq C_4\times C_4\),
its factorization into three prime discriminants
and eight possible candidates
for the metabelian Galois group \(G=\mathrm{Gal}(\mathrm{F}_2^\infty(K)/K)\)
of the \(2\)-class field tower of \(K\)
with length \(\ell_2(K)=2\)
in the notation of the SmallGroups library
\cite{BEO}
and ANUPQ package
\cite{GNO}.


\renewcommand{\arraystretch}{1.1}

\begin{table}[ht]
\caption{Discriminant of an imaginary quadratic field \(K\) with \(\ell_2(K)=2\)}
\label{tbl:UnexpImaginary2}
\begin{center}
\begin{tabular}{|r|c||c|c|}
\hline
 \(d_K\)     & Factorization            & \(\mathrm{Gal}(\mathrm{F}_2^2(K)/K)\) & \(\tau_4\) \\
\hline
 \(-2\,379\) & \((-3)\cdot 13\cdot 61\) & \(\langle 256,28\rangle-\#2;i-\#1;j\) & \((421)\)  \\
\hline
\end{tabular}
\end{center}
\end{table}


\noindent
The root paths of the occurring eight metabelian \(2\)-groups with order \(2048\) are given by
\eqref{eqn:Exp2RootPathUnexpI2},

\begin{equation}
\label{eqn:Exp2RootPathUnexpI2}
\begin{aligned}
\#1;j\to\#2;i\stackrel{2}{\to}\langle 256,28\rangle\stackrel{3}{\to}\langle 32,2\rangle, \quad i\in\lbrace 24,26\rbrace,\ 1\le j\le 4,
\end{aligned}
\end{equation}

\noindent
with respect to the lower exponent-\(2\) central series, and by
\eqref{eqn:RootPathUnexpI2},

\begin{equation}
\label{eqn:RootPathUnexpI2}
\begin{aligned}
\#1;j\to\#2;24\to\langle 512,791\rangle\stackrel{3}{\to}\langle 64,18\rangle\stackrel{2}{\to}\langle 16,2\rangle \\
\#1;j\to\#2;26\to\langle 512,792\rangle\stackrel{3}{\to}\langle 64,18\rangle\stackrel{2}{\to}\langle 16,2\rangle, \\
\end{aligned}
\end{equation}

\noindent
with respect to the (usual) lower central series.
The groups are of coclass \(6\).


\begin{theorem}
\label{thm:UnexpImaginary2} (Tower with length two.)
Let \(K\) be imaginary quadratic
with Artin pattern
\begin{equation}
\label{eqn:APUnexpImaginary2}
\begin{aligned}
\tau_0&=(22), \quad \tau_1=(311,311,411), \quad \varkappa_1=(J_0,K_2,K_3), \\
\tau_2&=(32,32,52,52,411,411;322), \quad \varkappa_2=(H_1,H_1,J_{31},J_{32},J_{21},J_{22};J_0), \quad \tau_4=(421).
\end{aligned}
\end{equation}
The \(2\)-class tower of \(K\) has length \(\ell_2(K)=2\)
and the metabelian Galois group \(G=\mathrm{Gal}(\mathrm{F}_2^\infty(K)/K)\)
is isomorphic to one of eight groups \(\langle 256,28\rangle-\#2;i-\#1;j\)
with order \(2048\) and relation rank \(d_2(G)=3\),
where \(i\in\lbrace 24,26\rbrace\) and \(1\le j\le 4\).
\end{theorem}


Here, we have \(\mathrm{Cl}_2(N_1)=(311)\) with order bigger than \(16\)
in the first layer already, whence \(\tau_4\) has \(2\)-rank at least three
\cite[Cor. 3.2]{BeSn},
independently of the \(2\)-class numbers \(\#\mathrm{Cl}_2(N_{1k})\).
Indeed, the precise structure is \(\tau_4=(421)\),
but nevertheless the tower length is only \(\ell_2(K)=2\).


\begin{example}
\label{exm:UnexpImaginary2}
The information in Formula
\eqref{eqn:APUnexpImaginary2}
occurs rather sparsely as Artin pattern
of imaginary quadratic fields \(K\) with \(2\)-class group of type \((4,4)\).
The absolutely smallest three discriminants \(d_K\) with this pattern are
\begin{equation}
\label{eqn:ExmUnexpImaginary2}
\begin{aligned}
-2\,379&=(-3)\cdot 13\cdot 61, \\
-9\,595&=5\cdot (-19)\cdot 101, \qquad
&-9\,955&=5\cdot (-11)\cdot 181.
\end{aligned}
\end{equation}
\end{example}


\begin{proof}
(Proof of Theorem
\ref{thm:UnexpImaginary2})
In the SmallGroups database
\cite{BEO}
the pattern \(\tau_0=(22)\) and \(\tau_1=(311,311,411)\), up to the first layer alone,
sets in with order \(512=2^9\).
It is hit by \(10\) groups
\(\langle 512,i\rangle\) with \(798\le i\le 801\), \(820\le i\le 822\) and \(851\le i\le 853\).
But none of the groups in the SmallGroups library hits the complete pattern up to the second layer,
including \(\tau_2=(32,32,52,52,411,411;322)\).
Some of the \(10\) candidates have too big components in their second layer pattern \(\tau_2\):
The three groups \(\langle 512,i\rangle\) with \(820\le i\le 822\) and common parent \(\langle 256,32\rangle\)
have inadequate Frattini invariants \((3211)\) instead of \((322)\),
the three groups \(\langle 512,i\rangle\) with \(851\le i\le 853\)and common parent \(\langle 256,38\rangle\)
have four inadequate invariants \((311)\) instead of \((32)\),
and they can be eliminated as parents for extending the pool by the \(p\)-group generation algorithm
\cite{Nm,Ob},
according to the \textit{monotony principle}
\cite[\S\ 1.7, Thm. 1.21, p. 79]{Ma2016b}.
The four groups \(\langle 512,i\rangle\) with \(798\le i\le 801\) and common parent \(\langle 256,28\rangle\)
possess compatible second layer \(\tau_2=(32,32,42,42,311,311;322)\)
but the descendants of the first two groups show stable components \((311)\) instead of \((411)\),
and thus the further construction of descendants can be stopped.

Instead, we select the group \(\langle 256,28\rangle\) as starting point for constructing an extended pool
of \(2\)-groups with orders bigger than \(512\), and thus outside of the SmallGroups database.
Descendants are generated with depth three (including great grand children)
and step sizes \(1\), \(2\) and \(3\) (including bi- and trifurcations to coclass bigger than \(5\)).
Among the \(4178\) groups generated in this way,
there are \(8\) hits of the complete pattern \(\tau_0=(22)\), \(\tau_1=(311,311,411)\),
\(\tau_2=(32,32,52,52,411,411;322)\) and \(\tau_4=(421)\) of abelian type invariants,
but without considering any capitulation patterns \(\varkappa_1\) and \(\varkappa_2\).
The eight groups are \(\langle 256,28\rangle-\#2;i-\#1;j\)
with order \(2048\) and admissible relation rank \(d_2=3\),
where \(i\in\lbrace 24,26\rbrace\) and \(1\le j\le 4\),
and they all share the desired capitulation type
\(\varkappa_1=(J_0,K_2,K_3)\) and \(\varkappa_2=(H_1,H_1,J_{31},J_{32},J_{21},J_{22};J_0)\).
Here, the trifurcation does not lead to non-metabelian groups with Artin pattern
\eqref{eqn:APUnexpImaginary2},
and the tower must necessarily be of length \(\ell_2(K)=2\),
since topologies of type parent --- child or fork --- sibling are discouraged.
\end{proof}


\section{Real quadratic fields of type \((4,4)\)}
\label{s:Real44}
\noindent
In \S\
\ref{s:Imaginary44},
we have seen that
the root path of all \(2\)-class field tower groups \(G=\mathrm{Gal}(\mathrm{F}_2^\infty(K)/K)\)
of imaginary quadratic fields \(K\) with \(\mathrm{Cl}_2(K)\simeq C_4\times C_4\)
contains at least one bifurcation,
and thus \(G\) must be at least of coclass \(\mathrm{cc}(G)=4\).
The deeper reason for this fact is a restriction
for the order of the capitulation kernel \(\ker(T_{K,N})\)
of \(K\) with respect to an unramified cyclic extension \(N/K\),
according to the Theorem on the Herbrand quotient of the unit group \(U_N\) of \(N\):
\begin{equation}
\label{eqn:Herbrand}
\#\ker(T_{K,N})=(U_K:\mathrm{N}_{N/K}(U_N))\cdot\lbrack N:K\rbrack.
\end{equation}
For an imaginary field \(K\),
the unit norm index is bounded by \(1\le (U_K:\mathrm{N}_{N/K}(U_N))\le 2\),
and Formula
\eqref{eqn:Herbrand}
yields the estimate \(4\le \#\ker(T_{K,N})\le 8\)
for an extension of degree \(\lbrack N:K\rbrack=4\),
whence \(\varkappa_2(K)\) cannot contain a total capitulation kernel.
However, for the finite \(2\)-groups \(G\)
near the root \(C_4\times C_4\),
\(\varkappa_2(G)\) always contains at least one total capitulation kernel,
whence the imaginary tower group \(G=\mathrm{Gal}(\mathrm{F}_2^\infty(K)/K)\) must be located rather far away from the root.

In the present section,
we shall see that
the root path of the \(2\)-class field tower group \(G=\mathrm{Gal}(\mathrm{F}_2^\infty(K)/K)\)
of a real quadratic field \(K\) with \(\mathrm{Cl}_2(K)\simeq C_4\times C_4\)
does not necessarily contain a bifurcation,
and thus \(G\) may be of the minimal possible coclass \(\mathrm{cc}(G)=3\).
For a real field \(K\) and an extension of degree \(\lbrack N:K\rbrack=4\), Formula
\eqref{eqn:Herbrand}
yields the less severe constraint \(4\le \#\ker(T_{K,N})\le 16\),
since \(1\le (U_K:\mathrm{N}_{N/K}(U_N))\le 4\),
due to the existence of the fundamental unit \(\eta\) in \(U_K\).
In fact, the real tower group \(G=\mathrm{Gal}(\mathrm{F}_2^\infty(K)/K)\)
may even coincide with the root:

\subsection{Abelian towers of Hilbert \(2\)-class fields}
\label{ss:RealLength1}
\noindent


\begin{theorem}
\label{thm:Abelian} (Single-stage tower.)
Let \(K\) be a real quadratic field
with Artin pattern
\begin{equation}
\label{eqn:APLength1}
\begin{aligned}
\tau_0&=(22), \quad \tau_1=(21,21,21), \quad \varkappa_1=(J_0,J_0,J_0), \\
\tau_2&=(2,2,2,2,2,2;11), \quad \varkappa_2=(G,G,G,G,G,G;G), \quad \tau_4=(0).
\end{aligned}
\end{equation}
Then the \(2\)-class field tower of \(K\) has length \(\ell_2(K)=1\)
and Galois group \(\mathrm{Gal}(\mathrm{F}_2^1(K)/K)\)
isomorphic to the abelian root \(\langle 16,2\rangle\simeq C_4\times C_4\)
with relation rank \(d_2=3\).
\end{theorem}


\begin{example}
\label{exm:Abelian}
The information in Formula
\eqref{eqn:APLength1}
cannot occur for imaginary quadratic fields but
seems to be a densely populated Artin pattern
of real quadratic fields \(K\) with \(2\)-class group of type \((4,4)\).
The smallest ten discriminants \(d_K\) with this pattern are
\begin{equation}
\label{eqn:ExmLength1}
\begin{aligned}
12\,104&=8\cdot 17\cdot 89, \qquad
&56\,648&=8\cdot 73\cdot 97, \\
83\,845&=5\cdot 41\cdot 409, \qquad
&136\,945&=5\cdot 61\cdot 449, \\
140\,488&=8\cdot 17\cdot 1033, \qquad
&173\,545&=5\cdot 61\cdot 569, \\
229\,445&=5\cdot 109\cdot 421, \qquad
&236\,249&=13\cdot 17\cdot 1069, \\
244\,936&=8\cdot 17\cdot 1801, \qquad
&276\,029&=13\cdot 17\cdot 1249.
\end{aligned}
\end{equation}
\end{example}


\begin{proof}
(Proof of Theorem
\ref{thm:Abelian})
According to item (1) of
\cite[Cor. 3.1, p. 1182]{BeSn},
the \(2\)-class field tower of \(K\) is abelian
with \(\mathrm{F}_2^\infty(K)=\mathrm{F}_2^1(K)\)
and \(\#\mathrm{Cl}_2(\mathrm{F}_2^1(K))=1\) (i.e., \(\tau_4=(0)\))
if and only if
\(\#\mathrm{Cl}_2(\mathrm{Fix}(H_i))=\frac{1}{2}\cdot\#\mathrm{Cl}_2(K)\),
for \(1\le i\le 3\)
(which is true for \(\tau_1=(21,21,21)\), since \(8=\frac{1}{2}\cdot 16\)).
Pattern
\eqref{eqn:APLength1}
is unique for the abelian group \(\langle 16,2\rangle\) in the SmallGroups library
\cite{BEO}.
\end{proof}


\subsection{Metabelian towers of Hilbert \(2\)-class fields}
\label{ss:RealLength2}
\noindent
In Table
\ref{tbl:Real2},
some discriminants \(0<d_K<10^6\) of real quadratic fields \(K\)
with \(\mathrm{Cl}_2(K)\simeq C_4\times C_4\) are shown together with
the factorization of \(d_K\) into prime discriminants
and the Galois group \(\mathrm{Gal}(\mathrm{F}_2^\infty(K)/K)\)
of the metabelian \(2\)-class field tower of \(K\)
in the notation of the SmallGroups library
\cite{BEO}
and ANUPQ package
\cite{GNO}.


\renewcommand{\arraystretch}{1.1}

\begin{table}[ht]
\caption{Discriminants of real quadratic fields \(K\) with \(\ell_2(K)=2\)}
\label{tbl:Real2}
\begin{center}
\begin{tabular}{|r|c||c|c|}
\hline
 \(d_K\)      & Factorization                & \(\mathrm{Gal}(\mathrm{F}_2^\infty(K)/K)\) & \(\tau_4\) \\
\hline
  \(58\,888\) & \(8\cdot 17\cdot 433\)       & \(\langle 32,4\rangle\)                    & \((1)\)    \\
  \(84\,972\) & \(4\cdot 3\cdot 73\cdot 97\) & \(\langle 64,20\rangle\)                   & \((2)\)    \\
 \(257\,045\) & \(5\cdot 101\cdot 509\)      & \(\langle 128,19\rangle\)                  & \((21)\)   \\
 \(232\,328\) & \(8\cdot 113\cdot 257\)      & \(\langle 256,281\rangle\)                 & \((31)\)   \\
  \(26\,245\) & \(5\cdot 29\cdot 181\)       & \(\langle 256,210\text{ or }211\rangle\)   & \((31)\)   \\
\hline
\end{tabular}
\end{center}
\end{table}

\noindent
The root paths of the occurring metabelian \(2\)-groups with orders up to \(512\) are given by
\eqref{eqn:RootPathR2}.

\begin{equation}
\label{eqn:RootPathR2}
\begin{aligned}
\langle 32,4\rangle\to\langle 16,2\rangle \\
\langle 64,20\rangle\to\langle 32,2\rangle\to\langle 16,2\rangle \\
\langle 128,19\rangle\to\langle 64,18\rangle\stackrel{2}{\to}\langle 16,2\rangle \\
\langle 256,281\rangle\to\langle 128,33\rangle\stackrel{2}{\to}\langle 32,2\rangle\to\langle 16,2\rangle
\end{aligned}
\end{equation}

\noindent
In accordance with the Shafarevich criterion,
all these groups have relation rank \(3\le d_2\le 4\).


\begin{theorem}
\label{thm:SporadicR} (Sporadic situation.)
Let \(K\) be a real quadratic field
with Artin pattern
\begin{equation}
\label{eqn:APSporadicR}
\begin{aligned}
\tau_0&=(22), \quad \tau_1=(31,31,22), \quad \varkappa_1=(K_1,K_1,K_1), \\
\tau_2&=(3,3,3,3,21,21;21), \quad \varkappa_2=(H_1,H_1,H_1,H_1,H_1,H_1;H_1), \quad \tau_4=(1).
\end{aligned}
\end{equation}
Then the \(2\)-class field tower of \(K\) has length \(\ell_2(K)=2\)
and metabelian Galois group \(\mathrm{Gal}(\mathrm{F}_2^\infty(K)/K)\)
isomorphic to the sporadic group \(\langle 32,4\rangle\)
with relation rank \(d_2=3\).
\end{theorem}


\begin{example}
\label{exm:SporadicR}
The information in Formula
\eqref{eqn:APSporadicR}
seems to be the most frequent Artin pattern
of real quadratic fields \(K\) with \(2\)-class group of type \((4,4)\).
The smallest ten discriminants \(d_K\) with this pattern are
\begin{equation}
\label{eqn:ExmSporadicR}
\begin{aligned}
58\,888&=8\cdot 17\cdot 433, \qquad
&73\,505&=5\cdot 61\cdot 241, \\
82\,433&=13\cdot 17\cdot 373, \qquad
&85\,969&=13\cdot 17\cdot 389, \\
110\,741&=37\cdot 41\cdot 73, \qquad
&116\,645&=5\cdot 41\cdot 569, \\
118\,001&=13\cdot 29\cdot 313, \qquad
&135\,505&=5\cdot 41\cdot 661, \\
143\,705&=5\cdot 41\cdot 701, \qquad
&149\,768&=8\cdot 97\cdot 193.
\end{aligned}
\end{equation}
\end{example}


\begin{proof}
(Proof of Theorem
\ref{thm:SporadicR})
According to item (2) of
\cite[Cor. 3.1, p. 1182]{BeSn},
the \(2\)-class field tower of \(K\) is metabelian
with \(\mathrm{F}_2^\infty(K)=\mathrm{F}_2^2(K)\)
and \(\#\mathrm{Cl}_2(\mathrm{F}_2^1(K))=2\) (i.e., \(\tau_4=(1)\))
if and only if
\(\#\mathrm{Cl}_2(\mathrm{Fix}(J_0))=\frac{1}{2}\cdot\#\mathrm{Cl}_2(K)\)
(which is true for \(\tau_2=(3,3,3,3,21,21;21)\), since \(8=\frac{1}{2}\cdot 16\)).
Pattern
\eqref{eqn:APSporadicR}
is unique for \(\langle 32,4\rangle\) in
\cite{BEO},
whereas \(\tau_1\) alone also appears for the child \(\langle 64,28\rangle\).
\end{proof}


\begin{theorem}
\label{thm:OtherSporadicR} (Another sporadic situation 1.)
Let \(K\) be a real quadratic field
with Artin pattern
\begin{equation}
\label{eqn:APOtherSporadicR}
\begin{aligned}
\tau_0&=(22), \quad \tau_1=(211,211,221), \quad \varkappa_1=(J_0,J_0,K_2), \\
\tau_2&=(21,21,21,21,22,22;211), \quad \varkappa_2=(G,G,G,G,H_2,H_2;H_2), \quad \tau_4=(2).
\end{aligned}
\end{equation}
Then the \(2\)-class field tower of \(K\) has length \(\ell_2(K)=2\)
and metabelian Galois group \(\mathrm{Gal}(\mathrm{F}_2^\infty(K)/K)\)
isomorphic to the sporadic group \(\langle 64,20\rangle\) with relation rank \(d_2=4\).
\end{theorem}


\begin{example}
\label{exm:OtherSporadicR}
The information in Formula
\eqref{eqn:APOtherSporadicR}
cannot occur for imaginary quadratic fields but
seems to be a sparsely populated Artin pattern
of real quadratic fields \(K\) with \(2\)-class group of type \((4,4)\).
The smallest four discriminants \(d_K\) with this pattern are
\begin{equation}
\label{eqn:ExmOtherSporadicR}
\begin{aligned}
84\,972&=4\cdot 3\cdot 73\cdot 97, \qquad
&177\,004&=4\cdot 17\cdot 19\cdot 137, \\
221\,001&=3\cdot 11\cdot 37\cdot 181, \qquad
&275\,356&=4\cdot 23\cdot 41\cdot 73.
\end{aligned}
\end{equation}
Note that all discriminants are divisible by four primes,
whence this situation with \(\tau_4=(2)\) lies outside of the scope of Theorems 6.2 and 6.3 in
\cite[pp. 1186--1191]{BeSn}
\end{example}


\begin{proof}
(Proof of Theorem
\ref{thm:OtherSporadicR})
According to item (3) of
\cite[Cor. 3.1, p. 1182]{BeSn},
the \(2\)-class field tower of \(K\) is metabelian
with \(\mathrm{F}_2^\infty(K)=\mathrm{F}_2^2(K)\)
and \(\mathrm{Cl}_2(\mathrm{F}_2^1(K))\) cyclic of order \(\ge 4\) (i.e., \(\tau_4=(e)\) with \(e\ge 2\))
if (but not only if)
\(\#\mathrm{Cl}_2(\mathrm{Fix}(H_i))=\#\mathrm{Cl}_2(K)\), for \(1\le i\le 2\),
\(\#\mathrm{Cl}_2(\mathrm{Fix}(H_3))\ge 2\cdot\#\mathrm{Cl}_2(K)\), and
\(\#\mathrm{Cl}_2(\mathrm{Fix}(J_{ij}))=\frac{1}{2}\cdot\#\mathrm{Cl}_2(K)\), for \(1\le i,j\le 2\)
(which is true for \(\tau_1=(211,211,221)\) and \(\tau_2=(21,21,21,21,22,22;211)\),
since \(32=2\cdot 16\) and \(8=\frac{1}{2}\cdot 16\)).
Pattern
\eqref{eqn:APOtherSporadicR}
is unique for \(\langle 64,20\rangle\) in the SmallGroups database
\cite{BEO},
whereas \(\tau_1\) alone also appears for orders \(128,256,512\).
\end{proof}


\begin{theorem}
\label{thm:OtherSporadic2R} (Another sporadic situation 2.)
Let \(K\) be a real quadratic field
with Artin pattern
\begin{equation}
\label{eqn:APOtherSporadic2R}
\begin{aligned}
\tau_0&=(22), \quad \tau_1=(211,211,311), \quad \varkappa_1=(J_0,J_0,K_1), \\
\tau_2&=(22,22,31,31,32,32;221), \quad \varkappa_2=(H_1,H_1,G,G,J_0,J_0;J_0), \quad \tau_4=(21).
\end{aligned}
\end{equation}
Then the \(2\)-class field tower of \(K\) has length \(\ell_2(K)=2\)
and metabelian Galois group \(\mathrm{Gal}(\mathrm{F}_2^\infty(K)/K)\)
isomorphic to the sporadic group \(\langle 128,19\rangle\) with relation rank \(d_2=3\).
\end{theorem}


\begin{example}
\label{exm:OtherSporadic2R}
The information in Formula
\eqref{eqn:APOtherSporadic2R}
seems to be a sparsely populated Artin pattern
of real quadratic fields \(K\) with \(2\)-class group of type \((4,4)\).
The smallest two discriminants \(d_K\) with this pattern are
\begin{equation}
\label{eqn:ExmOtherSporadic2R}
\begin{aligned}
257\,045&=5\cdot 101\cdot 509, \qquad
&259\,405&=5\cdot 29\cdot 1789.
\end{aligned}
\end{equation}
\end{example}


\begin{proof}
(Proof of Theorem
\ref{thm:OtherSporadic2R})
Pattern
\eqref{eqn:APOtherSporadic2R}
is unique for \(\langle 128,19\rangle\) in the SmallGroups library
\cite{BEO},
whereas \(\tau_1\) alone also appears for orders \(64,256,512\) and other groups of order \(128\).
Moreover, the transfer targets \(\tau_0\), \(\tau_1\), \(\tau_2\), \(\tau_4\)
also occur for \(\langle 128,18\rangle\),
which has, however, different tranfer kernels
\(\varkappa_1=(J_0,J_0,J_0)\) and \(\varkappa_2=(H_1,H_1,G,G,H_3,H_3;H_2)\).
\end{proof}


\begin{theorem}
\label{thm:OtherTreeR} (Coclass tree.)
Let \(K\) be a real quadratic field
with Artin pattern
\begin{equation}
\label{eqn:APOtherTreeR}
\begin{aligned}
\tau_0&=(22), \quad \tau_1=(211,311,311), \quad \varkappa_1=(J_0,K_2,J_0), \\
\tau_2&=(31,31,31,31,41,41;321), \quad \varkappa_2=(H_1,H_1,H_2,H_2,G,G;J_0), \quad \tau_4=(31).
\end{aligned}
\end{equation}
Then the \(2\)-class field tower of \(K\) has length \(\ell_2(K)=2\)
and metabelian Galois group \(\mathrm{Gal}(\mathrm{F}_2^\infty(K)/K)\)
isomorphic to the periodic group \(\langle 256,281\rangle\) with \(\zeta\simeq C_2\times C_2\) and \(d_2=3\)
on the tree \(\mathcal{T}^4(\langle 128,33\rangle)\).
\end{theorem}


\begin{example}
\label{exm:OtherTreeR}
The information in Formula
\eqref{eqn:APOtherTreeR}
seems to be a sparsely populated Artin pattern
of real quadratic fields \(K\) with \(2\)-class group of type \((4,4)\).
The smallest two discriminants \(d_K\) with this pattern are
\begin{equation}
\label{eqn:ExmOtherTreeR}
\begin{aligned}
232\,328&=8\cdot 113\cdot 257, \\
263\,993&=17\cdot 53\cdot 293 \quad 
\text{\cite[Exm. 7.3, pp. 1192]{BeSn}}.
\end{aligned}
\end{equation}
\end{example}


\begin{proof}
(Proof of Theorem
\ref{thm:OtherTreeR})
Searching the SmallGroups database
\cite{BEO}
for the pattern \(\tau_0=(22)\), \(\tau_1=(211,311,311)\) and \(\tau_2=(31,31,31,31,41,41;321)\)
yields \(4\) groups \(\langle 256,i\rangle\) with identifiers \(280\le i\le 283\) and
\(2\) groups \(\langle 512,j\rangle\) with identifiers \(1603\le j\le 1604\).
But only \(\langle 256,281\rangle\) possesses the correct capitulation type
\(\varkappa_1=(J_0,K_2,J_0)\)
and \(\varkappa_2=(H_1,H_1,H_2,H_2,G,G;J_0)\).
See Table
\ref{tbl:Contestants}.
Here, we have \(\mathrm{Cl}_2(N_1)=(211)\), \(\mathrm{Cl}_2(N_{1k})=(31)\), for \(k\in\lbrace 1,2\rbrace\),
and indeed \(\tau_4=(31)\) with \(2\)-rank two, according to
\cite[Cor. 3.2]{BeSn}.
\end{proof}


\renewcommand{\arraystretch}{1.1}

\begin{table}[ht]
\caption{Comparison of the contestants for Artin pattern \eqref{eqn:APOtherTreeR}}
\label{tbl:Contestants}
\begin{center}
{\small
\begin{tabular}{|c||c|c|c|c|}
\hline
 Vertex                      & \(\varkappa_1\) & \(\varkappa_2\)                 & \(\tau_4\) & \(\zeta\) \\
\hline
 \(\langle 256,280\rangle\)  & \(J_0,K_2,J_0\) & \(H_1,H_1,H_2,H_2,H_3,H_3;H_2\) & \((31)\)   & \((11)\) \\
 \(\langle 256,281\rangle\)  & \(J_0,K_2,J_0\) & \(H_1,H_1,H_2,H_2,G,G;J_0\)     & \((31)\)   & \((11)\) \\
 \(\langle 256,282\rangle\)  & \(J_0,K_2,J_0\) & \(H_1,H_1,H_2,H_2,H_2,H_2;J_0\) & \((31)\)   & \((11)\) \\
 \(\langle 256,283\rangle\)  & \(J_0,K_2,J_0\) & \(H_1,H_1,H_2,H_2,H_1,H_1;H_2\) & \((31)\)   & \((11)\) \\
\hline
 \(\langle 512,1603\rangle\) & \(J_0,K_2,J_0\) & \(H_1,H_1,H_2,H_2,H_3,H_3;H_2\) & \((41)\)   & \((1)\)  \\
 \(\langle 512,1604\rangle\) & \(J_0,K_2,J_0\) & \(H_1,H_1,H_2,H_2,H_3,H_3;H_2\) & \((41)\)   & \((1)\)  \\
\hline
\end{tabular}
}
\end{center}
\end{table}


\begin{theorem}
\label{thm:TreeVarR} (Another tree.)
Let \(K\) be a real quadratic field
with Artin pattern
\begin{equation}
\label{eqn:APTreeVarR}
\begin{aligned}
\tau_0&=(22), \quad \tau_1=(211,211,311), \quad \varkappa_1=(J_0,J_0,J_0), \\
\tau_2&=(22,22,31,31,32,32;321), \quad \varkappa_2=(H_1,H_1,G,G,H_3,H_3;H_2), \quad \tau_4=(31).
\end{aligned}
\end{equation}
Then the \(2\)-class field tower of \(K\) has length \(\ell_2(K)=2\)
and metabelian Galois group \(\mathrm{Gal}(\mathrm{F}_2^\infty(K)/K)\)
isomorphic to either \(\langle 256,210\rangle\) or \(\langle 256,211\rangle\)
with \(d_2=3\) on the tree \(\mathcal{T}^4(\langle 128,18\rangle)\).
\end{theorem}


\begin{example}
\label{exm:TreeVarR}
Only a single real quadratic fields \(K\) with \(2\)-class group of type \((4,4)\)
and Artin pattern in Formula
\eqref{eqn:APTreeVarR}
is known up to now.
The discriminant \(d_K\) of this field is
\(26\,245=5\cdot 29\cdot 181\).
\end{example}


\begin{proof}
(Proof of Theorem
\ref{thm:TreeVarR})
Searching the SmallGroups database
\cite{BEO}
for the pattern \(\tau_0=(22)\), \(\tau_1=(211,211,311)\) and \(\tau_2=(31,31,22,22,32,32;321)\)
yields \(4\) groups \(\langle 256,i\rangle\) with identifiers \(209\le i\le 212\).
But only \(\langle 256,210\rangle\) and \(\langle 256,211\rangle\) possess the correct capitulation type
\(\varkappa_1=(J_0,J_0,J_0)\)
and \(\varkappa_2=(H_1,H_1,G,G,H_3,H_3;H_2)\).
See Table
\ref{tbl:Contestants2}.
\end{proof}


\renewcommand{\arraystretch}{1.1}

\begin{table}[ht]
\caption{Comparison of the contestants for Artin pattern \eqref{eqn:APTreeVarR}}
\label{tbl:Contestants2}
\begin{center}
{\small
\begin{tabular}{|c||c|c|c|c|}
\hline
 Vertex                      & \(\varkappa_1\) & \(\varkappa_2\)             & \(\tau_4\) & \(\zeta\) \\
\hline
 \(\langle 256,209\rangle\)  & \(J_0,J_0,J_0\) & \(H_1,H_1,G,G,H_3,H_3;J_0\) & \((31)\)   & \((11)\) \\
 \(\langle 256,210\rangle\)  & \(J_0,J_0,J_0\) & \(H_1,H_1,G,G,H_3,H_3;H_2\) & \((31)\)   & \((11)\) \\
 \(\langle 256,211\rangle\)  & \(J_0,J_0,J_0\) & \(H_1,H_1,G,G,H_3,H_3;H_2\) & \((31)\)   & \((11)\) \\
 \(\langle 256,212\rangle\)  & \(J_0,J_0,J_0\) & \(H_1,H_1,G,G,H_3,H_3;J_0\) & \((31)\)   & \((11)\) \\
\hline
\end{tabular}
}
\end{center}
\end{table}

\subsection{Unexpected towers of Hilbert \(2\)-class fields with \(\ell_2(K)=2\)}
\label{ss:UnexpRealLength2}
\noindent 
Table
\ref{tbl:UnexpReal2}
shows a discriminant \(0<d_K<10^6\) of a real quadratic field \(K\)
with \(\mathrm{Cl}_2(K)\simeq C_4\times C_4\),
its factorization into prime discriminants
and the metabelian Galois group \(G=\mathrm{Gal}(\mathrm{F}_2^\infty(K)/K)\)
of the \(2\)-class field tower of \(K\) with unexpected length \(\ell_2(K)=2\)
in the notation of the SmallGroups library
\cite{BEO}
and ANUPQ package
\cite{GNO}.


\renewcommand{\arraystretch}{1.1}

\begin{table}[ht]
\caption{Discriminant of a real quadratic field \(K\) with \(\ell_2(K)=2\)}
\label{tbl:UnexpReal2}
\begin{center}
\begin{tabular}{|r|c||c|c|}
\hline
 \(d_K\)      & Factorization             & \(\mathrm{Gal}(\mathrm{F}_2^2(K)/K)\)         & \(\tau_4\) \\
\hline
  \(595\,561\) & \(17\cdot 53\cdot 661\)  & \(\langle 512,1465\rangle-\#1;2\) & \((222)\) \\
\hline
\end{tabular}
\end{center}
\end{table}

\noindent
The root path of the occurring metabelian \(2\)-group with order \(1024\) is given by
\eqref{eqn:UnexpRootPathR2}.

\begin{equation}
\label{eqn:UnexpRootPathR2}
\begin{aligned}
\#1;2\to\langle 512,1465\rangle\to\langle 256,237\rangle\to\langle 128,25\rangle\to\langle 64,19\rangle\stackrel{2}{\to}\langle 16,2\rangle \\
\end{aligned}
\end{equation}

\noindent
Here, the Shafarevich criterion cannot be used for a decision about the length,
since this metabelian group has relation rank \(d_2=3\),
and therefore can be \(2\)-tower group of real (and even imaginary) quadratic fields.


\begin{theorem}
\label{thm:UnexpTwoStagesR} (Two-stage tower.)
Let \(K\) be a real quadratic field
with Artin pattern
\begin{equation}
\label{eqn:APUnexpTwoStagesR}
\begin{aligned}
\tau_0&=(22), \quad \tau_1=(211,211,311), \quad \varkappa_1=(J_0,J_0,K_1), \\
\tau_2&=(31,31,31,31,311,422;2111), \quad \varkappa_2=(H_3,H_3,H_2,H_2,H_1,J_{11};H_1), \quad \tau_4=(222).
\end{aligned}
\end{equation}
Then the \(2\)-class tower of \(K\) has \(\ell_2(K)=2\)
and the metabelian Galois group \(G=\mathrm{Gal}(\mathrm{F}_2^\infty(K)/K)\)
is isomorphic to \(\langle 512,1465\rangle-\#1;2\) of order \(1024\).
\end{theorem}


\begin{example}
\label{exm:UnexpTwoStagesR}
\(595\,561\) is the single known discriminant \(0<d_K<10^6\)
of a real quadratic field \(K\) with \(\mathrm{Cl}_2(K)\simeq C_4\times C_4\)
and Artin pattern in Formula
\eqref{eqn:APUnexpTwoStagesR}.
It is mentioned in
\cite[Exm. 7.4, p. 1192]{BeSn}
without explicit invariant \(\tau_4=(222)\).
\end{example}


\begin{proof}
(Proof of Theorem
\ref{thm:UnexpTwoStagesR})
Searching the SmallGroups database
\cite{BEO}
for the pattern \(\tau_0=(22)\) and \(\tau_1=(211,211,311)\), up to the first layer alone,
yields \(2\) groups of order \(64\),
\(13\) groups of order \(128\), for instance \(\langle 128,25\rangle\),
\(21\) groups of order \(256\), for instance \(\langle 256,237\rangle\), and
\(56\) groups of order \(512\), for instance \(\langle 512,1465\rangle\).
But none of the groups in the SmallGroups library hits the complete pattern up to the second layer,
including \(\tau_2=(31,31,31,31,311,422;2111)\),
because \(\tau_2(6)\) takes inadequate values in \(\lbrace 311,321,322\rbrace\),
which are too small but compatible with \((422)\), by the monotony principle.

We can select the group \(\langle 128,25\rangle\),
resp. its child \(\langle 256,237\rangle\),
resp. its grand child \(\langle 512,1465\rangle\),
as the starting point for selectively constructing an extended pool
of \(2\)-groups with orders bigger than \(512\), and thus outside of the SmallGroups database.
We generate descendants with
depth three (including great grand children),
resp. depth two (including grand children),
resp. depth one (only immediate children),
and step size \(1\) alone,
since there are no bifurcations to coclass bigger than \(4\).
Among the \(54\), resp. \(47\), resp. \(4\), groups generated in this manner,
there are always exactly \(2\) hits of the complete pattern \(\tau_0=(22)\), \(\tau_1=(211,211,311)\),
\(\tau_2=(31,31,31,31,311,422;2111)\) and \(\tau_4=(222)\) of abelian type invariants,
but without considering any capitulation patterns \(\varkappa_1\) and \(\varkappa_2\).

Table
\ref{tbl:PolUnexpReal2}
shows the development of the Artin pattern along the root path
of the four descendants \(\langle 512,1465\rangle-\#1;i\) with \(1\le i\le 4\).


\renewcommand{\arraystretch}{1.1}

\begin{table}[ht]
\caption{Development of the Artin pattern along the root path}
\label{tbl:PolUnexpReal2}
\begin{center}
{\small
\begin{tabular}{|c||c|c|c|c|c|}
\hline
 Vertex                      & \(\tau_1\)       & \(\varkappa_1\) & \(\tau_2\)                   & \(\varkappa_2\)                    & \(\tau_4\) \\
\hline
 \(\langle 16,2\rangle\)     & \(21,21,21\)     & \(J_0,J_0,J_0\) & \(2,2,2,2,2,2;11\)           & \(G,G,G,G,G,G;G\)                  & \((0)\)    \\
 \(\langle 64,19\rangle\)    & \(211,211,211\)  & \(J_0,J_0,J_0\) & \(31,31,31,31,31,31;211\)    & \(H_3,H_3,H_2,H_2,H_1,H_1;G\)      & \((2)\)    \\
 \(\langle 128,25\rangle\)   & \(211,211,311\)  & \(J_0,J_0,K_1\) & \(31,31,31,31,311,311;2111\) & \(H_3,H_3,H_2,H_2,H_1,H_1;H_1\)    & \((21)\)   \\
 \(\langle 256,237\rangle\)  & \(211,211,311\)  & \(J_0,J_0,K_1\) & \(31,31,31,31,311,321;2111\) & \(H_3,H_3,H_2,H_2,H_1,H_1;H_1\)    & \((211)\)  \\
 \(\langle 512,1465\rangle\) & \(211,211,311\)  & \(J_0,J_0,K_1\) & \(31,31,31,31,311,322;2111\) & \(H_3,H_3,H_2,H_2,H_1,H_1;H_1\)    & \((221)\)  \\
 \(\#1;1\)                   & \(211,211,311\)  & \(J_0,J_0,K_1\) & \(31,31,31,31,311,422;2111\) & \(H_3,H_3,H_2,H_2,H_1,J_{12};H_1\) & \((222)\)  \\
 \(\#1;2\)                   & \(211,211,311\)  & \(J_0,J_0,K_1\) & \(31,31,31,31,311,422;2111\) & \(H_3,H_3,H_2,H_2,H_1,J_{11};H_1\) & \((222)\)  \\
 \(\#1;3\)                   & \(211,211,311\)  & \(J_0,J_0,K_1\) & \(31,31,31,31,311,322;2111\) & \(H_3,H_3,H_2,H_2,H_1,H_1;H_1\)    & \((222)\)  \\
 \(\#1;4\)                   & \(211,211,311\)  & \(J_0,J_0,K_1\) & \(31,31,31,31,311,322;2111\) & \(H_3,H_3,H_2,H_2,H_1,H_1;H_1\)    & \((222)\)  \\
\hline
\end{tabular}
}
\end{center}
\end{table}


The groups \(\langle 512,1465\rangle-\#1;i\) with \(3\le i\le 4\) can be eliminated,
because \(\tau_2(6)=(322)\) instead of \((422)\).
Now the capitulation type
\(\varkappa_1=(J_0,J_0,K_1)\) and \(\varkappa_2=(H_3,H_3,H_2,H_2,H_1,J_{11};H_1)\)
is taken into account,
and we are unambiguously led to \(\langle 512,1465\rangle-\#1;2\),
since \(\langle 512,1465\rangle-\#1;1\) is discouraged.
\end{proof}


\subsection{Smallest three-stage tower of Hilbert \(2\)-class fields}
\label{ss:SmallestLength3}
\noindent
To find a \(2\)-class field tower \(\mathrm{F}_2^\infty(K)/K\) of length \(\ell_2(K)=3\)
with \textit{minimal degree} over a quadratic field \(K\) having \(\mathrm{Cl}_2(K)\simeq C_4\times C_4\)
we changed the technique of investigation.
Instead of beginning with Artin patterns of fields \(K\)
we started with a search for non-metabelian \(2\)-groups \(G\) with \(G/G^\prime\simeq C_4\times C_4\)
of minimal order in the SmallGroups database
\cite{BEO}.
The result was a batch of \(16\) groups \(G\) of order \(256=2^8\)
with identifiers \(\langle 256,i\rangle\), where \(469\le i\le 484\).
It turned out that all these groups share a common
second derived quotient (metabelianization) \(G/G^{\prime\prime}\),
and that the latter coincides with their common parent \(\pi(G)\simeq\langle 128,122\rangle\).

This situation was called a \textit{simple tree topology} of type \textit{parent and child} in
\cite[\S\ 5, pp. 89--92]{Ma2016b}.
The root path of the occurring \(2\)-groups with orders up to \(256\) is given by
\eqref{eqn:SmallestLength3}.
It is free from any bifurcations.
All descendants are of step size \(1\), and the coclass \(3\) remains stable.

\begin{equation}
\label{eqn:SmallestLength3}
\begin{aligned}
\langle 256,i\rangle\to\langle 128,122\rangle\to\langle 64,23\rangle\to\langle 32,2\rangle\to\langle 16,2\rangle
\end{aligned}
\end{equation}

\noindent
Since \(d_2(\pi(G))=5\), the tower length in the following theorem is strictly proven.


\begin{theorem}
\label{thm:SmallestLength3} (Minimal three-stage tower.)
Let \(K\) be a real quadratic field
with Artin pattern
\begin{equation}
\label{eqn:APSmallestLength3}
\begin{aligned}
\tau_0&=(22), \quad \tau_1=(211,211,211), \quad \varkappa_1=(J_0,J_0,J_0), \\
\tau_2&=(21,21,21,21,211,2111;1111), \quad \varkappa_2=(G,G,G,G,G,H_3;H_1), \quad \tau_4=(111).
\end{aligned}
\end{equation}
The \(2\)-class tower of \(K\) has length \(\ell_2(K)=3\)
and the metabelianization \(G/G^{\prime\prime}=\mathrm{Gal}(\mathrm{F}_2^2(K)/K)\)
of the Galois group \(G=\mathrm{Gal}(\mathrm{F}_2^\infty(K)/K)\)
is isomorphic to \(\langle 128,122\rangle\) with coclass \(3\) and relation rank \(5\).
The group \(G\) itself could be one of the twelve immediate descendants
\(\langle 256,i\rangle\) with relation rank \(4\) and coclass \(3\), where
\(470\le i\le 472\) or \(474\le i\le 476\) or \(478\le i\le 480\) or \(482\le i\le 484\).
\end{theorem}


\begin{example}
\label{exm:SmallestLength3}
In reverse direction, we used the common Artin pattern \(\mathrm{AP}(G)\)
of all these \(12\) non-metabelian groups \(G\) and of their metabelianization \(G/G^{\prime\prime}\)
in Formula
\eqref{eqn:APSmallestLength3}
for finding quadratic fields \(K\) having this pattern \(\mathrm{AP}(K)=\mathrm{AP}(G)\).
Indeed, we discovered two suitable \textit{real} quadratic fields with the following minimal discriminants \(d_K\):
\begin{equation}
\label{eqn:ExmSmallestLength3}
\begin{aligned}
150\,997&=(-7)\cdot (-11)\cdot 37\cdot 53, \\
224\,652&=(-4)\cdot (-3)\cdot 97\cdot 193.
\end{aligned}
\end{equation}
\end{example}


\begin{proof}
(Proof of Theorem
\ref{thm:SmallestLength3})
Searching the SmallGroups database
\cite{BEO}
for the pattern \(\tau_0=(22)\) and \(\tau_1=(211,211,211)\), up to the first layer alone,
yields \(1\) group of order \(32\),
\(5\) groups of order \(64\),
\(12\) groups of order \(128\),
\(52\) groups of order \(256\), and
\(218\) groups of order \(512\).
The complete pattern up to the second layer,
including \(\tau_2=(21,21,21,21,211,2111;1111)\),
is hit by
a single group of order \(128\), namely \(\langle 128,122\rangle\),
\(16\) groups of order \(256\), namely \(\langle 256,i\rangle\) with \(469\le i\le 484\),
\(24\) groups of order \(512\), namely \(\langle 512,j\rangle\) with \(1967\le j\le 1990\), and
\(16\) groups of order \(1024\).
They all share the correct capitulation type
\(\varkappa_1=(J_0,J_0,J_0)\) and \(\varkappa_2=(G,G,G,G,G,H_3;H_1)\),
and minimal possible coclass \(3\).
The group of order \(128\) is metabelian with relation rank \(5\),
all the other groups have derived length \(3\).
Only \(12\), resp. \(23\), groups of orders \(256\) and \(1024\), resp. \(512\),
are terminal with mandatory relation rank \(4\).
\end{proof}


\section{Conclusion}
\label{s:Summary}
\noindent
We have shown impressively that it is sufficient to compute
the logarithmic abelian type invariants of \(2\)-class groups \(\mathrm{Cl}_2(N)\)
and the capitulation kernels
of unramified abelian extensions \(N/K\)
in the first layer (\(\tau_1\) and \(\varkappa_1\))
and second layer (\(\tau_2\) and \(\varkappa_2\))
over the quadratic field \(K\),
that is, with absolute degrees \(4\) and \(8\).
Together with the foregiven base layer \(\tau_0=(22)\),
corresponding to the general assumption that \(\mathrm{Cl}_2(K)\simeq C_4\times C_4\),
this information was the input for the pattern recognition process
in our group theoretic search and the preamble of all theorems.

We never had to determine data for the third layer (\(\tau_3\) and \(\varkappa_3\)),
let alone for the Hilbert \(2\)-class field in the fourth layer (\(\tau_4\)),
with absolute degrees \(16\) and \(32\).
(Observe that the capitulation \(\varkappa_4\) is total,
according to the principal ideal theorem by Artin and Furtw\"angler.)

We emphasize that, on the contrary, we obtained the
logarithmic abelian type invariants \(\tau_4\)
of the \(2\)-class group \(\mathrm{Cl}_2(\mathrm{F}_2^1(K))\)
of the Hilbert \(2\)-class field of \(K\)
as the structure of the commutator subgroup \(G^\prime\)
of the sifted metabelian \(2\)-group \(G=\mathrm{Gal}(\mathrm{F}_2^2(K)/K)\).
This is a significant progress in the
determination of deeper arithmetical information
compared with the paper of Benjamin and Snyder
\cite{BeSn},
where Pari/GP was unable to compute \(\tau_4=(222)\) in Example 7.4 on page 1192.

The length of the \(2\)-class field tower of \(K\)
turned out to be either \(\ell_2(K)=3\) or unexpectedly also \(\ell_2(K)=2\) (in Theorem
\ref{thm:UnexpImaginary2}
and Theorem
\ref{thm:UnexpTwoStagesR}),
when the rank of the \(2\)-class group \(\mathrm{Cl}_2(\mathrm{F}_2^1(K))\) is three.
We only found metabelian towers with length \(\ell_2(K)=2\),
when the rank of the \(2\)-class group \(\mathrm{Cl}_2(\mathrm{F}_2^1(K))\) is two,
although \(\ell_2(K)=3\) should be possible, according to Blackburn
\cite[p. 1164]{BeSn}.


\section{Outlook on unsolved problems}
\label{s:Outlook}
\noindent
We do not want to pretend the impression that this paper
has covered all aspects of the \(2\)-class field tower
of quadratic fields \(K\) with \(2\)-class group \(\mathrm{Cl}_2(K)\simeq C_4\times C_4\).
We restricted our investigations mainly to the cases where Corollary 3.2
\cite[p. 1182]{BeSn}
by Benjamin and Snyder gives
a warranty for a \(2\)-class group \(\mathrm{Cl}_2(\mathrm{F}_2^1(K))\) of rank two
of the Hilbert \(2\)-class field \(\mathrm{F}_2^1(K)\) of \(K\).
In these cases, we always found a two-stage tower.
We also touched upon some towers with three or even more stages,
making clear that these questions lead to
the difficult problem of higher abelian quotient invariants.
We intentionally avoided the study of Artin patterns with high complexity, as given in Table
\ref{tbl:HighComplexity},
because there it is probably impossible to determine merely
the metabelianization \(G/G^{\prime\prime}\) of the \(2\)-tower group \(G\).


\renewcommand{\arraystretch}{1.1}

\begin{table}[ht]
\caption{Artin patterns with high complexity}
\label{tbl:HighComplexity}
\begin{center}
{\small
\begin{tabular}{|c||c|c|c|}
\hline
 Discriminant & \(\tau_1\)      & \(\tau_2\)                            & \(\tau_4\)    \\
\hline
  \(-6\,328\) & \(311,411,611\) & \(42,42,72,72,411,411;532\)           & \((622)\)     \\
  \(-7\,672\) & \(311,311,311\) & \(311,311,321,321,422,422;21111\)     & \((22222)\)   \\
 \(-10\,803\) & \(311,311,311\) & \(321,321,321,321,322,322;21111\)     & \((3321111)\) \\
 \(-12\,595\) & \(221,221,711\) & \(221,221,721,822,2211,3221;61111\)   & \((72221)\)   \\
 \(-19\,947\) & \(221,221,411\) & \(221,221,421,522,2211,4411;31111\)   & \((43321)\)   \\
 \(-20\,155\) & \(221,221,411\) & \(221,221,421,432,2211,3311;31111\)   & \((43332)\)   \\
 \(-27\,955\) & \(221,221,411\) & \(522,522,2211,2211,2221,4311;31111\) & \((5432211)\) \\
\hline
\end{tabular}
}
\end{center}
\end{table}


\section{Acknowledgements}
\label{s:Thanks}

\noindent
We gratefully acknowledge that our research was supported by the Austrian Science Fund (FWF):
project P 26008-N25.
We are indebted to Professor M. F. Newman from the Australian National University in Canberra
for drawing our attention to the possibility of a four-stage tower in Theorem
\ref{thm:ThreeStagesI2}
and the finitely presented \(2\)-groups in Theorem
\ref{thm:Limits}.



\end{document}